\documentclass[onefignum,onetabnum]{siamart190516}

\usepackage{subfig}
\usepackage{amsmath}
\usepackage{amssymb}
\usepackage{lipsum}
\usepackage{amsfonts}
\usepackage{graphicx}
\usepackage{epstopdf}
\usepackage{algorithmic}
\ifpdf
  \DeclareGraphicsExtensions{.eps,.pdf,.png,.jpg}
\else
  \DeclareGraphicsExtensions{.eps}
\fi

\usepackage{enumitem}

\newsiamremark{remark}{Remark}
\newsiamremark{hypothesis}{Hypothesis}
\crefname{hypothesis}{Hypothesis}{Hypotheses}
\newsiamthm{claim}{Claim}

\theoremstyle{plain}
\theoremheaderfont{\normalfont\sc}
\theorembodyfont{\normalfont\itshape}
\theoremseparator{.}
\newtheorem*{thm}{Standing Assumption}

\headers{Analysis of a Networked Bivirus Epidemic Model}{M. Ye, B. D. O. Anderson, and J. Liu}

\title{Convergence and Equilibria Analysis of a Networked Bivirus Epidemic Model\thanks{Submitted to the editors: \today.
\funding{M. Ye was supported in part by the Western Australian Government through the Premier's Science Fellowship Program. B.D.O. Anderson was supported in part by the Australian Research Council (ARC) under grant DP160104500 and DP190100887.}}}

\author{Mengbin Ye\thanks{School of Electrical Engineering, Computing and Mathematical Sciences, Curtin University, Perth 6102, WA, Australia 
  (\email{mengbin.ye@curtin.edu.au}).}
\and Brian D.O. Anderson\thanks{School Engineering, Australian National University, Acton 2601, A.C.T., Australia
  (\email{brian.anderson@anu.edu.au}).}
\and Ji Liu\thanks{Department of Electical and Computer Engineering, State University of New York (SUNY) at Stony Brook, NY, United States of America.
  (\email{ji.liu@stonybrook.edu}).}}

\usepackage{amsopn}

\ifpdf
\hypersetup{
  pdftitle={Analysis of a Networked Bivirus Epidemic Model},
  pdfauthor={M. Ye, B. D. O. Anderson, and J. Liu}
}
\fi

\begin{document}

\maketitle

% REQUIRED
\begin{abstract}
This paper studies a networked bivirus model, in which two competing viruses spread across a network of interconnected populations; each node represents a population with a large number of individuals. The viruses may spread through possibly different network structures, and an individual cannot be simultaneously infected with both viruses. Focusing on convergence and equilibria analysis, a number of new results are provided. First, we show that for networks with generic system parameters, there exist a finite number of equilibria. Exploiting monotone systems theory, we further prove that for bivirus networks with generic system parameters, then convergence to an equilibrium occurs for all initial conditions, except possibly for a set of measure zero. Given the network structure of one virus, a method is presented to construct an infinite family of network structures for the other virus that results in an infinite number of equilibria in which both viruses coexist. Necessary and sufficient conditions are derived for the local stability/instability of boundary equilibria, in which one virus is present and the other is extinct. A sufficient condition for a boundary equilibrium to be almost globally stable is presented. Then, we show how to use monotone systems theory to generate conclusions on the ordering of stable and unstable equilibria, and in some instances identify the number of equilibria via rapid simulation testing. Last, we provide an analytical method for computing equilibria in networks with only two nodes, and show that it is possible for a bivirus network to have an unstable coexistence equilibrium and two locally stable boundary equilibria.
%   {\color{red} To be rewritten} Here is an example SIAM \LaTeX\ article. This can be used as a
%   template for new articles.  Abstracts must be able to stand alone
%   and so cannot contain citations to the paper's references,
%   equations, etc.  An abstract must consist of a single paragraph and
%   be concise. Because of online formatting, abstracts must appear as
%   plain as possible. Any equations should be inline. 250 word limit.
\end{abstract}

% REQUIRED

% REQUIRED (see https://mathscinet.ams.org/mathscinet/msc/msc2020.html to check)

\section{Introduction}
Mathematical models of epidemic spreading processes have been of interest to the broad scientific community for decades~\cite{review}, and have recently come into the limelight as a result of the ongoing COVID-19 pandemic~\cite{Vespignani2020model}. In the context of infectious diseases, such models are of interest to predict the dynamics of the disease and the course of an outbreak. One may seek to draw conclusions on whether the disease will eventually disappear or become endemic, examine the impact of key characteristics such as infection and recovery rates in shaping the epidemic, and design control strategies to reduce or stop the spread~\cite{nowzari2016epidemics,pare2017epidemic,review,mei2017epidemics_review,ye2021_PH_TAC}.

Among the different modeling frameworks, compartmental models have become especially popular; different compartments indicate particular health states and each individual in a large population may move between compartments due to the infectious nature of the disease. A classical compartmental model is the Susceptible--Infected--Susceptible (SIS) model, in which each individual is either healthy and susceptible (S) to the disease, or infected (I) but can recover with no immunity~\cite{review}. If immunity is permanently acquired after recovery, then the Susceptible--Infected--Removed (SIR) framework is used~\cite{review}. For a particular compartmental framework, both stochastic and deterministic models exist, and are often related by mean-field approximation~\cite{nowzari2016epidemics}. Although not as realistic in capturing the probabilistic nature of infectious disease transmissions, deterministic models remain popular as they offer a balance in terms of analytical tractability, simulation implementation and modeling accuracy. 

More recently, attention has grown on models that study the spread of two or more diseases/viruses within the same population (the literature most commonly refers to these as bi- or multi-virus models). Depending on the application of interest, the viruses may be assumed to compete against, reinforce, or weaken one another~\cite{sahneh2014competitive,newman2013interacting,granell2013dynamical}. Competitive virus models assume that if an individual is infected by one virus, then they cannot be infected by any other virus, and are especially popular. Different variants of SIS- or SIR-like models also have broad applications beyond infectious diseases, including competing ideas, decisions, and internet memes~\cite{sahneh2014competitive,wei2013competing,karrer2011competing}.

This paper focuses on the ``bivirus model'' in the SIS framework, which considers two competing viruses, called virus 1 and virus 2, spreading across a network of interconnected populations on possibly different network structures~\cite{liu2019bivirus,prakash2012winner,sahneh2014competitive,yang2018bivirus,santos2015bivirus}. Such a scenario may reflect two strains of a disease spreading in a population, such as gonorrhea and a strain of partially drug-resistant gonorrhea~\cite{carlos1,carlos2}.
% {\color{cyan}A real-world example consists of two distinct sexually transmitted pathogens or two competing strains of the same pathogen (e.g., gonorrhea and a strain of partially drug-resistant gonorrhea) spreading in a heterosexually active population \cite{carlos1,carlos2}. It is worth mentioning that in \cite{carlos1} the whole population is assumed to be ``genetically and behaviorally homogeneous except for the gender of individuals in the population'' and thus it splits into male and female populations, which can be viewed as a real-world situation of a low order system considered in Subsection~\ref{sssec:n2}.}
To discuss existing works, we first introduce the reproduction number of virus 1 and virus 2 as $\mathcal{R}^1$ and $\mathcal{R}^2$, respectively. 
In general, the infection and recovery rate parameters are assumed to be heterogeneous among the nodes of the network, and then $\mathcal{R}^i$ for $i \in \{1, 2\}$ is a nontrivial but computable function of the recovery and infection rates. For homogeneous parameters, $\mathcal{R}^i$ simplifies to the ratio of the infection and recovery rate for virus $i$, multiplied by the spectral radius of the adjacency matrix of the virus $i$ spreading network.
% {\color{black}When the infection and recovery rate parameters are homogeneous among the nodes of the network, $\mathcal{R}^i$ is the ratio of the infection and recovery rate for virus $i \in \{1, 2\}$, multiplied by the spectral radius of the weighting matrix associated to the network structure for the spread of virus $i$. When the parameters are heterogeneous among the nodes, $\mathcal{R}^i$ is a nontrivial but computable function of the recovery and infection rates and the weighting matrix.} 
In the case that only virus $i$, for $i\in \{1,2\}$, is present in the population, $\mathcal{R}^i \leq 1$ or $\mathcal{R}^i > 1$ determines if virus $i$ asymptotically disappears from or persists in  the population, respectively. Much of the analysis in existing works focuses on the assumption that $\mathcal{R}^i \leq 1$ for at least one of $i=1$ or $i = 2$; in this case, either one or both viruses will be eliminated asymptotically~\cite{liu2019bivirus,yang2018bivirus,santos2015bivirus}. 
The works of \cite{sahneh2014competitive,santos2015bivirus,santos2015bivirus_conference} consider $\mathcal{R}^i > 1$ for both $i=1,2$. Assuming homogeneous infection and recovery rates among individuals for each virus, but allowing generic network structures, \cite{sahneh2014competitive} deals only with local stability of equilibria. In contrast,  \cite{santos2015bivirus,santos2015bivirus_conference} assume general infection rates but homogeneous recovery rates, and focus on global stability of specific equilibria of interest.
% {\color{black}The works of \cite{sahneh2014competitive,santos2015bivirus,santos2015bivirus_conference} consider $\mathcal{R}^i > 1$ for both $i=1,2$. Assuming homogeneous infection and recovery rates among individuals for each virus, but allowing generic network structures, \cite{sahneh2014competitive} analyses the equilibria and their local stability properties (including the emergence of \textit{coexistence} equilibria in which both virus are present). In \cite{santos2015bivirus}, a sufficient condition on general infection rates (with retention of homogeneous recovery rates) is provided which guarantees that for all initial conditions, a specific virus survives in the network while the other is eliminated. A different sufficient condition is given in~\cite{santos2015bivirus_conference}.}

This paper considers the deterministic bivirus model, with the two viruses spreading on possibly different network structures across the same set of nodes (populations). We allow heterogeneous infection and recovery parameters, which yields new dynamical phenomena, such as the possibility of an unstable coexistence equilibrium, where both viruses are present. We establish several novel results on the system dynamics and associated equilibria.
\begin{enumerate}
    \item We formally prove that a bivirus system with generic values of infection and recovery parameters has a finite number of equilibria. Prior to our work, the finiteness of equilibria has only been proved for special parameter values~\cite{liu2019bivirus,janson2020networked}, or assumed to be true with no rigorous arguments~\cite{sahneh2014competitive}. We also prove that one can conduct equilibria analysis by assuming unity recovery rates without loss of generality, thus reducing the parameter space dimension. The work of \cite{santos2015bivirus_conference,santos2015bivirus} made a similar assumption but without justification.
    \item We use a {\color{black}coordinate} transformation to demonstrate that the bivirus system has monotone trajectories, simplifying the proof in \cite{santos2015bivirus}. Going beyond \cite{santos2015bivirus}, we explicitly connect the bivirus dynamics to the monotone systems literature~\cite{hirsch1985systems,smith1988monotone_survey,smith2008monotone_book}. We significantly extend the existing results by providing a general conclusion on the limiting dynamical behaviour of the bivirus networked system: a generic bivirus system converges to an equilibrium among the finite number of equilibria for all initial conditions except possibly a set of measure zero. Thus, no chaos is possible, and limit cycles (if they exist) are nonattractive. Our result differs from~\cite{santos2015bivirus,santos2015bivirus_conference,liu2019bivirus}, which establish sufficient conditions for a particular equilibrium to be (almost) globally attractive. By further exploiting the literature, we establish several conclusions on the ordering of stable and unstable equilibria.
    % This is the first general convergence result in the case that $\mathcal{R}^i > 1, i =1,2$.
    \item We identify a class of system parameters that yields a connected set containing an infinite number of equilibria, comprising an interval of a straight line, in which both viruses coexist, covering a much broader set of parameters than \cite{liu2019bivirus}. We term such a connected set a ``line of equilibria'' for brevity, and show that this line of equilibria is locally exponentially attractive. Then, simple necessary and sufficient conditions are given for checking whether boundary equilibria, in which one virus is present and the other is absent, are locally stable or unstable. A sufficient condition for one of the boundary equilibria to be almost globally attractive is provided.
    % An existing result establishes that for specific special parameter sets, {\color{black}there exists a connected set containing an infinite number of equilibria, comprising an interval of a straight line, in} which both viruses coexist~\cite{liu2019bivirus}. {\color{black}We term such a connected set a ``line of equilibria'' for brevity.} We expand this by identifying a broader class of system parameters that yields a line of equilibria, and we further show that this line of equilibria is locally exponentially attractive. Then, simple necessary and sufficient conditions are given for checking whether boundary equilibria, in which one virus is present and the other is absent, are {\color{black}locally} stable or unstable.
    \item For networks with only two nodes, an analytic method to compute coexistence equilibria. We then provide numerical examples showcasing the different stability and instability configurations for the boundary and coexistence equilibria. We report a system which has an unstable coexistence equilibrium, and two locally stable boundary equilibria with nontrivial regions of attraction, highlighting the nontrivial equilibria and convergence properties of networked bivirus systems.
    % ; all examples of coexistence equilibria reported in the literature thus far have been locally stable.
    % \item Exploiting the monotone systems property, we show that essentially all trajectories of the bivirus system can be bounded by two trajectories with special initial conditions. This enables rapid simulation testing and a concrete condition for networks with only two nodes on when a single equilibrium exists in which both virus are present.
\end{enumerate}

% The outline is not required, but we show an example here.
The paper is organized as follows. \Cref{sec:prelims} details the mathematical preliminaries and presents the bivirus model. The main results are detailed in \cref{sec:general}, and conclusions are presented in \cref{sec:conclusions}.

\section{Preliminaries}\label{sec:prelims}

In this section, we introduce notation, and specify the nature of the system under consideration. The set-up is drawn from that in \cite{liu2019bivirus}, and we also review the key conclusions of \cite{liu2019bivirus,sahneh2014competitive,santos2015bivirus,santos2015bivirus_conference}. This section also includes a review of some key results in monotone systems theory \cite{smith2008monotone_book,smith1988monotone_survey}, a tool we shall use in this paper. 

\subsection{Notation and related}\label{ssec:linear_algebra}
Given a natural number $n$, define the index set $[n] = \{1, 2, \hdots, n\}$. For real vectors $x,y\in\mathbb R^n$, with entries $x_i, y_i , i\in [n]$, we write $x\geq y$ iff $x_i\geq y_i\;\forall i$, we write $x>y$ iff $x_i\geq y_i\;\forall i$ and $x\neq y$, and we write $x\gg y$ iff $x_i>y_i\;\forall i$. For matrices $A,B$ of the same dimension, the notation $A\geq B, A>B, A\gg B$ mean the same thing as the corresponding inequalities relating ${\rm{vec}}(A), {\rm{vec}}(B)$; we say that $A$ is nonnegative if $A \geq 0_{n\times n}$. The $n$-vectors of all 1's and 0's are denoted by ${\bf 1}_n, {\bf 0}_n$ respectively, and the $n$-dimensional identity matrix is denoted $I_n$. The sets $\mathbb R^n_{\geq 0}$ and $\mathbb R^n_{>0}$ denote $\{x\in\mathbb R^n: x\geq {\bf 0}_n\}$ and $\{x\in\mathbb R^n: x\gg {\bf 0}_n\}$ respectively.  The set $\{x\in\mathbb R^n,{\bf 0}_n\leq x\leq {\bf 1}_n\}$ will be denoted by $\Xi_n$. 

Suppose $A$ is a square matrix. Then $\rho(A)$ and $s(A)$ will denote the spectral radius of $A$ and the greatest real part of any eigenvalue of $A$, respectively. It is Hurwitz if $s(A)<0$. We say that $A$ is reducible iff there is a permutation matrix $P$ such that $P^{\top}AP$ is block upper triangular; otherwise $A$ is said to be irreducible.
When a nonnegative $A$ is irreducible, and $Ax=y$ for $x>{\bf 0}_n, y>{\bf 0}_n$, $y$ cannot have a zero entry in every position where $x$ has a zero entry; equivalently, $A$ is reducible if there exist $x>{\bf 0}_n,y>{\bf 0}_n$ with the same set of zero entries. 
A square matrix $A$ is a Metzler matrix if all off-diagonal entries are nonnegative. For an irreducible Metzler $A$, and by an extension of the Perron-Frobenius theorem \cite{horn1994topics_matrix}, $s(A)$ is a simple eigenvalue and the only eigenvalue with this real part,  and the corresponding eigenvector can be taken to have all positive entries, while the eigenvectors corresponding to other eigenvalues do not have this property. A square matrix $A$ is an $M$-matrix if $-A$ is Metzler and all eigenvalues of $A$ have positive real parts except for any at the origin; we say $A$ is a singular or nonsingular $M$-matrix if it has an eigenvalue at the origin or if its eigenvalues have strictly positive real parts, respectively~\cite{horn1994topics_matrix}. Further key properties, detailed in \cite[Theorem ~2.1]{varga2009matrix_book} and \cite[Theorem 2.3 and Theorem 4.6]{berman1979nonnegative_matrices}, are:
\begin{enumerate}
    \item 
For Metzler $A$, $s(A)$ depends monotonically on any entry of $A$. 
    \item
 For a matrix $F$ that is a (singular) $M$-matrix, and any positive diagonal $D$, $DF$ is also a (singular) $M$-matrix.
 \item For an irreducible nonnegative matrix $B$ and positive diagonal matrix $D$, there holds i) $s(-D+B) > 0 \Leftrightarrow \rho(D^{-1}B) > 1$, ii) $s(-D+B) = 0 \Leftrightarrow \rho(D^{-1}B) = 1$ and iii) $s(-D+B) < 0 \Leftrightarrow \rho(D^{-1}B) < 1$.
 \end{enumerate}

We shall work with weighted directed graphs $\mathcal G=(\mathcal {V,E},A)$, with $\mathcal V=\{1,2,\dots,n\}$ as the vertex set, $\mathcal E \subseteq \mathcal V \times \mathcal V$ the edge set and $A$ a nonnegative $n\times n$ square matrix. Further, $a_{ij}>0$ if and only if $(j,i) \in \mathcal E$, which connotes the existence of a directed edge from node $j$ to node $i$. Strong connectivity of $\mathcal G$ is the property that there exists a path connecting any two nodes, and is equivalent to $A$ being irreducible~\cite{berman1979nonnegative_matrices}.  

\subsection{The bivirus equations}\label{ssec:bivirus}
To begin, we review the classical Susceptible--Infected--Susceptible (SIS) network model~\cite{lajmanovich1976SISnetwork}. We consider $n>1$ populations of individuals, with each population being of large and constant size. Each individual has two possible health states, being healthy but susceptible (S), or infected (I) with the virus. Individuals who recover are assumed to do so with no immunity, being immediately susceptible again to infection. Defining the fraction $x_i\in[0,1]$ of population $i \in [n]$ who are infected, the SIS model posits that
\begin{equation}\label{eq:onevirus}
    \dot x_i(t)=-\delta_ix_i(t)+(1-x_i(t))\sum_{j=1}^n\beta_{ij}x_j(t),\quad i \in [n]
\end{equation}
where $\delta_i > 0$ is the recovery rate of population $i$ and $\beta_{ij}\geq 0$ is the transmission rate of the virus from infected individuals in population~$j$ to susceptible individuals in population~$i$. With $x=[x_1,x_2,\dots,x_n]^{\top}$, $X={\rm{diag}}(x)$, $B=(\beta_{ij})$, and \mbox{$D={\rm{diag}}(\delta_1,\delta_2,\dots,\delta_n)$}, the system equation becomes
\begin{equation}\label{eq:sis}
    \dot x(t)=[-D+\left (I_n-X(t)\right)B]x(t)
\end{equation}
The forward invariance of the set $\Xi_n$ for \eqref{eq:sis} is well known~\cite{lajmanovich1976SISnetwork}, and this guarantees that $x_i$, for all $i\in[n]$, retains its physical meaning. The graph $\mathcal{G} = (\mathcal{V}, \mathcal{E}, B)$ captures the network structure over which the virus spreads, and the standard assumption that $B$ is irreducible is equivalent to strong connectivity of $\mathcal{G}$; the virus can reach any population $i$ from any other population $j$ through a path on $\mathcal{G}$.

Given the irreducibility of $B$, a complete convergence characteristic can be determined by the reproduction number $\mathcal{R} \triangleq \rho(D^{-1}B)$. If $\mathcal{R} \leq 1$, which is equivalent to $s(-D+B)\leq 0$, the only equilibrium of \eqref{eq:sis} is the healthy equilibrium ${\bf 0}_n$ and it is globally asymptotically stable for $x(0)\in\Xi_n$ (exponentially stable if $s(-D+B)<0$). If $\mathcal{R} > 1$, i.e. $s(-D+B)>0$, there is in addition to the healthy equilibrium (which is unstable) a unique nonzero/endemic\footnote{We call such an equilibrium endemic because the virus is present in at least one node.} equilibrium $\bar x\gg {\bf 0}_n$ which is globally attractive for $x(0)\in \Xi_n\setminus {\bf{0}}_n$ \cite{lajmanovich1976SISnetwork,mei2017epidemics_review}. The equation \eqref{eq:sis} is the mean-field approximation of a stochastic discrete-time system with $2n$ states. Detailed discussion of the relationship between \eqref{eq:sis} and the stochastic model can be found in~\cite{nowzari2016epidemics,vanMeighem2009_virus,van2015accuracy}, and we do not explore this aspect further.

This SIS framework can be extended to study the scenario in which there are circulating two viruses, termed virus 1 and virus 2, as in~\cite{liu2019bivirus,sahneh2014competitive,santos2015bivirus}. They are competing, in that an individual infected with one virus has immunity to an infection from the other virus, but like the single virus framework, an infected individual can recover with no immunity and immediately become susceptible again to infection from either virus. Similar to \eqref{eq:sis}, the bivirus dynamics presented below is a mean-field approximation of a $3^n$-state Markov chain model. Due to space limitations, we do not discuss the details of this approximation. The reader is referred to \cite[Section~II]{liu2019bivirus} and also the comprehensive Chapter~5 in the PhD thesis~\cite{santos2014PhDThesis} for two different treatments.
% The virus dynamics are captured as a continuous-time deterministic system, and can be considered as a mean-field approximation of a stochastic discrete-time system; we do not discuss this relationship further, as it is detailed in \cite{liu2019bivirus}, and focus instead on the analysis of the system.   

Neighbor relationships for each virus are modeled by a directed graph
% \footnote{The graph is usually different for the two viruses, to avoid a nongeneric situation.}
between the populations (nodes), corresponding to the graph vertices, with an edge from node $j$ to node $i$ denoting the direction in which virus transmission can occur; the nonnegative infection rates $\beta^1_{ij}, \beta ^2_{ij}$ for $i,j\in[n]$ capture the transmission rates for virus 1 and virus 2, respectively. Each group $i \in [n]$ also has associated with it positive $\delta^1_i, \delta^2_i$ corresponding to the rate of recovery from virus 1 and virus 2, respectively. 

With virus 1 and virus 2, we distinguish the vectors of fractions of infected individuals as $x^1(t)\in\mathbb R^n$ and $x^2(t)\in\mathbb R^n$, respectively. The corresponding system equations become
\begin{subequations}\label{eq:twovirus}
\begin{align}
\dot x^1(t)=[-D^1+\left(I_n-X^1(t)-X^2(t)\right)B^1]x^1(t), \label{eq:virus_1} \\
\dot x^2(t)=[-D^2+\left(I_n-X^1(t)-X^2(t)\right)B^2]x^2(t). \label{eq:virus_2}
\end{align}
\end{subequations}

% {\color{cyan}The above model is an approximated model resulting from a mean-field approximation on
% a $3^n$-state Markov chain model. Detailed derivation of the mean-field approximation can be found in \cite{liu2019bivirus} and Chapter~5 in \cite{santos2014PhDThesis}}.

In order for $x^1, x^2$ to have physical meaning, observe that in addition to the requirement that ${\bf{0}}_n\leq x^i(t)\leq {\bf{1}}_n$ for $i=1,2$, we require that $x^1(t)+x^2(t)\leq {\bf{1}}_n$, the quantity on the left corresponding to the vector of fractions of individuals infected by either virus. Following a similar proof to the single virus case, an invariant set for the bivirus model can be identified which satisfies such requirements.
\begin{lemma}[{\cite[Lemma 8]{liu2019bivirus}}]\label{lem:looseinvariant}
With the above notation, suppose that the initial conditions for \eqref{eq:twovirus} satisfy ${\bf{0}}_n\leq x^i(0)\leq {\bf{1}}_n$ for $i=1,2$, and $x^1(0)+x^2(0)\leq {\bf{1}}_n$. Then for all $t>0$, there holds ${\bf{0}}_n\leq x^i(t)\leq {\bf{1}}_n$ for $i=1,2$ and $x^1(t)+x^2(t)\leq {\bf{1}}_n$.
\end{lemma}

% There are of course many results in \cite{liu2019bivirus}; some are for general systems, and some are for very specialized systems, where the assumed infection and healing rates are nongeneric. It is not simply the case that the assumed rates are nongeneric, with conclusions relating to stability for the nongeneric rates perhaps pointing to what happens with generic rates, but a number of the conclusions are themselves nongeneric, and thus \textit{not} indicative of what happens with generic rates. Our interest is more in the latter and we now highlight several conclusions for generic systems of \cite{liu2019bivirus}. 

% Despite the number of results in \cite{liu2019bivirus,sahneh2014competitive,yang2018bivirus,santos2015bivirus}, several major open questions are left unanswered concerning equilibria and convergence, which we will address in this paper.
In order to discuss existing results in the literature such as \cite{liu2019bivirus,sahneh2014competitive,yang2018bivirus,santos2015bivirus}, and gain an appreciation of some open questions that remain unanswered, we introduce a key assumption to hold in the rest of the paper that parallels one typically imposed when analyzing the single virus case\footnote{Relaxation of this assumption can of course be considered separately.}.
% In the absence of the assumption, there comes the possibility that certain groups, initially uninfected, could remain uninfected due to there being no influence, or path of infection, from other groups. Such a situation could of course be examined separately. 

\begin{thm}\label{standing_assump}
% {\color{red} Numbering system for standing assump to be cleaned up by Ben I think.} 
The matrices $B^1,B^2$ are nonnegative and irreducible. The matrices $D^1, D^2$ are positive diagonal. 
\end{thm}  

Positive diagonal $D^1, D^2$ implies that for any population (node) $i \in \mathcal{V}$, the recovery rate against each of the two viruses is strictly positive, but may differ between the viruses and between populations. By associating $B^i$ with the graph $\mathcal{G}^i = (\mathcal{V}, \mathcal E^i, B^i)$, for $i = 1,2$, one can view $\mathcal{G}^1$ and $\mathcal{G}^2$ as the graphs capturing the network over which virus 1 and virus 2 spread, respectively. The assumption on irreducibility is equivalent to both $\mathcal{G}^1$ and $\mathcal{G}^2$ being strongly connected, but we allow the edges and associated infection rates to differ between the two graphs. 

We first explain some reasoning behind several existing results of \cite{liu2019bivirus,santos2015bivirus,yang2018bivirus}, as opposed to reproducing the proofs. Consider action of a single virus, say virus 1, and the effect on $x^1$ in two cases: when virus 2 is not present, and when virus 2 is present. First, let us define the reproduction numbers of virus 1 and virus 2 as $\mathcal{R}^1 \triangleq \rho((D^1)^{-1}B^1)$ and $\mathcal{R}^2 \triangleq \rho((D^2)^{-1}B^2)$, respectively. Examination of \eqref{eq:onevirus} and \eqref{eq:virus_1} shows that the presence of $x^2 \geq {\bf 0}_n$ serves to slow down the increase in $x^1$ due to infection, 
% i.e. to reduce an additive term contributing to $\dot x^1(t)$. Put another way, 
because there necessarily holds $(I_n-X^1(t)-X^2(t))B^1x^1(t)\leq (I_n-X^1(t))B^1x^1(t)$. This means that if in the absence of virus 2, the healthy state $x^1 = {\bf 0}_n$ is an attractive equilibrium state for virus 1, i.e. any nonzero $x^1(0)$ yields $x^1(t) \to {\bf 0}_n$, then $x^1 = {\bf 0}_n$ must \textit{a fortiori} remain an attractive equilibrium in the presence of virus 2. This implies then that if $\mathcal{R}^1 \leq 1$, any equilibrium $(\bar x^1,\bar x^2)$ of the bivirus system \eqref{eq:twovirus} will necessarily be of the form $({\bf 0}_n, \bar x^2)$, and consequently in this case $\bar x^2$ is necessarily an equilibrium of the single virus system applicable to virus 2 alone due to \eqref{eq:virus_2}. And then, there will further hold $\bar x^2={\bf 0}_n$ if and only if $\mathcal{R}^2 \leq 1$. (Of course, one can interchange the roles of the two viruses). These observations effectively cover \cite[Theorems 1 to 3]{liu2019bivirus}, \cite{yang2018bivirus}, and part of \cite{santos2015bivirus}. The remaining theorems in \cite{liu2019bivirus} focus on equilibria analysis of systems \eqref{eq:twovirus} with nongeneric parameter matrices $D^i, B^i, i = 1,2$. As identified in \cite{liu2019bivirus} and as we will further explore in the sequel, stability and equilibria of nongeneric parameters are \textit{not always} indicative of what happens for generic parameters.  

The astute reader will recognize that the observations above also clarify that the condition $\mathcal{R}^1 > 1$ \textit{and} $\mathcal{R}^2 > 1$ remains one of great interest, since this is \textit{not} covered by a collapse to  two single virus problems.
% {\color{black}Indeed, this paper will often focus on advancing results this condition; no general convergence results exist, with only local stability results for homogeneous parameters in \cite{sahneh2014competitive} and limited analysis of the equilibria in \cite{liu2019bivirus}.}
Note that when both these conditions hold, there will still necessarily be three equilibria $({\bf 0}_n, {\bf 0}_n), (\bar x^1, {\bf 0}_n)$ and $({\bf 0}_n,\bar x^2)$; we term the first the \textit{healthy} equilibrium and the latter two \textit{boundary} equilibria of the bivirus system. We remark that $\bar x^1 \gg {\bf 0}_n$ and $\bar x^2 \gg {\bf 0}_n$ are separately the unique endemic equilibria for each of the two single virus systems. This is evident since for $x^2(0)={\bf 0}_n$, the underlying equations \eqref{eq:twovirus} imply $x^2(t)={\bf 0}_n$ for all $t$, and the evolution of $x^1(t)$ from some $x^1(0)$ is the same as what would occur with the corresponding single virus system in \eqref{eq:onevirus}. This also implies that $(\bar x^1, {\bf 0}_n)$ is the only equilibrium where virus 1 is present and virus 2 is not present. The same holds correspondingly with virus 1 and 2 interchanged. 

The works \cite{sahneh2014competitive,santos2015bivirus,santos2015bivirus_conference} consider $\mathcal{R}^1 > 1$ \textit{and} $\mathcal{R}^2 > 1$, and allow $\mathcal{G}^1$ and $\mathcal G^2$ to differ. In \cite{sahneh2014competitive}, infection and recovery parameters are assumed homogeneous among the nodes for any one virus, but can be different between the two viruses. The existence and local stability of equilibria is studied using linearization and perturbation methods, coupled with extensive numerical simulations~\cite{sahneh2014competitive}. In \cite{santos2015bivirus} and \cite{santos2015bivirus_conference}, differing sufficient conditions on infection parameters but with homogeneous recovery parameters are presented which guarantee the boundary equilibrium $({\bf 0}_n,\bar x^2)$ is attractive for all initial conditions satisfying $x^2(0) > {\bf 0}_n$ (and thus a similar condition exists for $({\bf 0}_n,\bar x^2)$ to be attractive). However, a number of important open questions remain for generic bivirus systems. For instance, i) are there a finite number of equilibria? and ii) can a dichotomy be established for the typical limiting behaviour? In this paper, we answer these questions and more, extensively and rigorously expanding on existing results.

% Apart from $({\bf 0}_n, {\bf 0}_n)$ being unstable, little is known about the (in)stability of other equilibria of \eqref{eq:twovirus}  or equivalently, $s(-D^1+B^1)>0$ and $s(-D^2+B^2)>0$.

\subsection{Monotone systems theory} 

As background for the results in later sections, we outline the concept of monotone systems \cite{smith1988monotone_survey, smith2008monotone_book} and recall some key results. Indeed, in the sequel, we use a standard method to demonstrate that the bivirus system is a monotone system, offering an alternative and simpler proof than the approach in \cite{santos2015bivirus}.  The discussion in this subsection however is in terms of a general system 
\begin{equation}\label{eq:sys}
    \dot x=F(x),\quad x\in \mathbb{R}^n
\end{equation}
on a convex, % {\footnote{Relaxation of convexity to `order convexity' is possible as explained in \cite[Remark~3]{smith1988monotone_survey}.}}
open set $U\subset \mathbb R^n$. We assume henceforth that conditions on $F$ exist which guarantee existence and uniqueness of solutions for all time, and indeed the Jacobian $dF_x$ exists at every point on allowed trajectories. Let $\phi_t(x_0)$ denote the solution $x(t)$ of \eqref{eq:sys} at time $t$ when $x(0)=x_0\in U$.  There is special interest in {\color{black}the} behavior of the \textit{difference} between trajectory pairs when confined to particular orthants of $\mathbb R^n$, and we begin with the central definition. 

Let $m = (m_1, m_2, \hdots, m_n)$, with $m_i\in \{0,1\}$ for all $i\in [n]$, be a prescribed sequence, and associate with it an orthant
\begin{equation}
    K_m=\{x\in\mathbb R^n:(-1)^{m_i}x_i\geq 0,\forall i\in[n]\}.
\end{equation}
We say $x\leq_{K_m} y$ and $x\ll_{K_m}y$ if $y-x\in K_m$ and $y-x\in{\rm{int}}(K_m)$ respectively. Here, $\rm{int}(\cdot)$ denotes the interior of a set.

Then, \eqref{eq:sys} is termed a  {\bf{Type $K_m$ monotone system}} if whenever $x_0,y_0\in U$ satisfy $x_0\leq_{K_m} y_0$, there holds $\phi_t(x_0)\leq_{K_m} \phi_t(y_o)$. {\color{black}That is}, $\phi_t$ preserves the partial ordering $\leq _{K_m}$ for $t\geq 0$.

There is a straightforward condition for the monotone property, in terms of $dF_x$, see~\cite[Lemma~2.1]{smith1988monotone_survey} and~\cite[p. 424]{hirsch1985systems}: 

\begin{lemma}\label{lem:kamke}[Kamke--M\"uller Condition]
With the above notation including conditions on $F$ and $U$, $\phi_t(x_0)$ preserves the partial ordering $\leq_{K_m}$ for $t\geq 0$ if and only if, with definition $P_m={\rm{diag}}\left((-1)^{m_1},(-1)^{m_2},\dots,(-1)^{m_n}\right)$, the matrix $P_mdF_xP_m$ is a Metzler matrix for all $x\in U$.
% , where 
% \[P_m={\rm{diag}}\left((-1)^{m_1},(-1)^{m_2},\dots,(-1)^{m_n}\right).
% \]
If in addition $P_mdF_xP_m$ is irreducible for all $x\in U$, then $x_0\leq_{K_m}y_0$ implies $\phi_t(x_0)\ll_{K_m}\phi_t(y_0)$ for all $t>0$.
\end{lemma}

Given this matrix characterization, it makes sense to define an {\bf{irreducible monotone system}} as one for which $P_mdF_xP_m$ is irreducible for all $x\in U$. 

One of the most important results concerning irreducible monotone systems is that only certain forms of limiting behavior are permitted, provided that there are known to be only a finite number of equilibria, see \cite[\color{black}Theorems~2.5 and 2.6]{smith1988monotone_survey}:

\begin{lemma}\label{lem:monotone_convergence}
Suppose $\dot x=F(x)$ is an irreducible type $K_m$ monotone system with $\mathcal M$ an open, bounded and positively invariant set. Suppose the closure of $\mathcal M$, denoted $\overline{\mathcal M}$, contains a finite number of equilibria $x^i$, with basins of attraction $\mathcal B(x^i)$. {\color{black}Let the set of all equilibria in $\overline{\mathcal M}$ by denoted by $E$. Then the set $\mathcal Y\in \mathcal M$ of points $x$ for which $\phi_t(x)$ does not converge to an equilibrium is of Lebesgue measure zero, and $\cup_{x^i\in E}{\rm{int}}(\mathcal B(x^i))\cap \overline{\mathcal M}$ is open and dense in $\mathcal M$. }
% \[
% \cup_{x^i}{\rm{int}}(\mathcal B(x^i))\cap \overline{\mathcal M}
% \]
% is open and dense in $\mathcal M$.
% , and any point outside this set lies on a nonattractive limit cycle. If there is a unique equilibrium $x^*\in\mathcal M$ (and no equilibrium in $\overline{\mathcal M}\setminus \mathcal M$), then $\phi_t(x_0)$ converges to $x^*$  from all $x(0)\in\mathcal M$. {\color{red} For discussion with Ben. Can we perhaps dispose of the requirement to state "no equilibrium in $\overline{\mathcal M}\setminus \mathcal M$" if we impose some kind of requirement on trajectories at the boundary? It seems to be a nuisance assumption, especially for the bivirus problem where we may have an equilibrium on the boundary defined by a single virus equilibrium of the form $\bar x^1,0)$ say. I can't even remember why this assumption came in. Do we know of an example where it is critical?}
\end{lemma}
The above result establishes that the monotone system cannot exhibit chaos, and any limit cycles in  $\mathcal{M}$ must be nonattractive. In other words, convergence to an equilibrium occurs for all initial conditions in {\color{black}$\mathcal{M}$} except possibly a set $\mathcal{Y}$ of measure zero. Note that if $x^i$ is unstable, then ${\rm{int}}(\mathcal B(x^i)) = \emptyset$. This implies that convergence will occur to a \textit{stable} equilibrium, except possibly for a set $\mathcal{Z}$ of initial conditions of measure zero which would yield convergence to an unstable equilibrium if one exists (a saddle point or a source). Clearly $\mathcal{Z} \subset \mathcal{Y}$, if in fact the two subsets are not empty.
% , and Note that when there is only one equilibrium, even more is known, see \cite[Proposition~4]{ye2019_PH_submit}.

\section{Main results}\label{sec:general}

In this section, we will present a series of new findings on the bivirus system. In order to keep the focus on the new results and to aid exposition and discussion, all proofs are presented in the Appendix. Moreover, and unless explicitly stated otherwise, the \nameCref{standing_assump} is assumed to hold throughout.

In much of this section we will appeal to the equilibrium equations associated with \eqref{eq:twovirus}. With an overbar denoting an equilibrium, these equations are
\begin{subequations}\label{eq:equil}
\begin{align}
[-D^1+(I_n-\bar X^1-\bar X^2)B^1]\bar x^1&=0, \label{eq:equil_1}  \\
[-D^2+(I_n-\bar X^1-\bar X^2)B^2]\bar x^2&=0. \label{eq:equil_2}
\end{align}
\end{subequations}
Since the outcomes when $x^i(0) ={\bf 0}_n$ for some $i$ are fully understood from the single virus dynamics, we focus on initial conditions in the set
\begin{equation*}
    \Delta = \{(x^1, x^2) \;\vert\; {\bf 0}_n < x^i < {\bf 1}_n \text{ for } i = 1,2, \text{and } x^1+x^2 \leq {\bf 1}_n \}.
\end{equation*}
Thus, when we refer to a globally attractive equilibrium, it is with respect to initial conditions in $\Delta$, and an almost globally attractive equilibrium excludes initial conditions in a set of measure zero that is a subset of $\Delta$.

\subsection{General result on equilibria and convergence}
In this subsection, we establish general properties of the equilibria of \eqref{eq:twovirus}, and provide a general convergence result. First, \eqref{eq:equil} can be used to establish the following.
\begin{lemma}\label{lem:interiorequilibria}
%  Under the \Cref{standing_assump}, 
 Any solution $(\bar x^1,\bar x^2)$ of the equilibrium equations \eqref{eq:equil} with $\bar x^i\geq {\bf 0}_n$ is such that for each $i = 1,2$, either $\bar x^i={\bf 0}_n$ or $\bar x^i\gg {\bf 0}_n$. Moreover, $\bar x^1+\bar x^2\ll{\bf{1}}_n$. Suppose $(\bar x^1,\bar x^2)$ and $(\bar x^1,\tilde x^2)$ are two equilibria. If $\bar x^1\gg {\bf 0}_n$, then $\bar x^2=\tilde x^2$. If $\bar x^1 = {\bf 0}_n$, and $\bar x^2 \gg {\bf 0}_n$ and $\tilde x^2 \gg {\bf 0}_n$, then $\bar x^2=\tilde x^2$. 
\end{lemma}

We term any equilibrium $(\bar x^1,\bar x^2)$ with $\bar x^1\gg {\bf 0}_n$ and $\bar x^2\gg {\bf 0}_n$ a \textit{coexistence} equilibrium.
This lemma, whose proof is given in \cref{app:pf_interiorequib}, restricts substantially the equilibria which can lie on the boundary of the set $\{x^1,x^2 \in \mathbb{R}^n_{\geq 0} | x^1+x^2 \leq {\bf 1}_n\}$. There are corresponding restrictions on the trajectories, as expressed in the following, which strengthens \Cref{lem:looseinvariant}, and is essentially obtained as a result of taking into account the \nameCref{standing_assump} (see \cref{app:pf_strictinvariant} for the proof).

\begin{lemma}\label{lem:strictinvariant}
Suppose that the bivirus equation set \eqref{eq:twovirus} has initial conditions satisfying $(x^1(0), x^2(0)) \in \Delta$.
% Suppose further \Cref{standing_assump} is in force.
Then for all finite\footnote{The results of \Cref{lem:strictinvariant,lem:interiorequilibria} do not preclude the possibility of $\lim_{t\to\infty} x^i(t) = {\bf 0}_n$. Thus, it is possible to converge to an equilibrium in which one or both viruses are not present, but at any finite time $t > 0$, either the virus is present in all nodes, or present in none. Such observations are consistent with those of the single virus SIS model~\cite{lajmanovich1976SISnetwork}.} $t>0$, there holds ${\bf{0}}_n\ll x^i(t)\ll {\bf{1}}_n$ for $i=1,2$ and $x^1(t)+x^2(t)\ll {\bf{1}}_n$.
\end{lemma}

It is well known that if $\bar x$ is an equilibrium of \eqref{eq:sys}, then the (local) stability of $\bar x$ can often be determined through examination of the eigenvalues of the Jacobian $dF_x(\bar x)$. If the linear system $\dot z=dF_x(\bar x)z$ is (exponentially) asymptotically stable or is unstable, the same is true of \eqref{eq:sys}, at least locally around $\bar x$~\cite[Section~5.8]{sastry1999nonlinearbook}. 

In the case of the bivirus system, the Jacobian provides information not only concerning equilibria. Denoting by $J$ the Jacobian of the bivirus system, it is straightforward to verify that, with $\tilde B^i(\bar x^i)={\rm{diag}}(B^i\bar x^i)$,
\small
\begin{align}\label{eq:jacobian}
    J(\bar x^1,\bar x^2) & =
    \begin{bmatrix}
    -D^1+(I_n-\bar X^1-\bar X^2)B^1-\tilde B^1(\bar x^1)&-\tilde B^1(\bar x^1)\\-\tilde B^2(\bar x^2)&-D^2+(I_n-\bar X^1-\bar X^2)B^2-\tilde B^2(\bar x^2)
    \end{bmatrix}
\end{align}
\normalsize

The Jacobian at a general point on a trajectory is the same as that in \eqref{eq:jacobian}, save that $\bar x^i, \bar X^i$ are replaced by $x^i,X^i$. Let $P \triangleq {\rm{diag}}(I_n,-I_n)$. Then observe (now with $\tilde B^i={\rm{diag}}(B^ix^i)$) that

\begin{equation}\label{eq:jacob_transform}
\small{
PJ(x^1,x^2)P=\begin{bmatrix}
    -D^1+(I_n-X^1-X^2)B^1&0\\0&-D^2+(I_n-X^1- X^2)B^2
    \end{bmatrix}
    +
    \begin{bmatrix}
    -\tilde B^1&\tilde B^1\\\tilde B^2&-\tilde B^2
    \end{bmatrix}
    }
\end{equation}
{\color{black}It is immediately clear that this matrix is a Metzler matrix}. However, more is true, as the following result illustrates (the proof is given in \cref{app:pf_jacobirred}):

\begin{lemma}\label{lem:jacobianirreducible}
Suppose $(x^1(0), x^2(0))\in \Delta$. Then the matrix $PJ(x^1(t),x^2(t))P$  evaluated at an arbitrary point along the trajectories of \eqref{eq:twovirus} with $t<\infty$, and the matrix $PJ(\bar x^1,\bar x^2)P$ for an equilibrium satisfying $\bar x^1\gg {\bf 0}_n$,  $\bar x^2\gg {\bf 0}_n$, are irreducible.
\end{lemma}

Thus the bivirus system given the \nameCref{standing_assump} is an irreducible monotone system of type $K_m$ where $m=({\bf 0}_n^\top, {\bf 1}_n^\top)$. While monotonicity of the bivirus trajectories with homogeneous recovery rates was established in \cite[Theorem 18]{santos2015bivirus}, the proof via \cref{lem:jacobianirreducible} is significantly more direct. We also explicitly make a connection to the monotone systems literature\footnote{Monotonicity of the single virus system \eqref{eq:onevirus} was identified several decades ago, see~\cite{hirsch1982monotonesystems_I,smith1988monotone_survey}.} and competitive species models~\cite{hsu1996competitive,lam2016remark,carlos2}. Indeed, the bivirus system is a \textit{cooperative system} in the terminology of \cite{smith1988monotone_survey,hirsch1985systems} (not to be confused with competitive species models). A key contribution of this paper is to show how to leverage the monotone systems literature, including the seminal contributions by Morris M.~Hirsh in the 1980s, to establish a number of new results for networked bivirus systems. Our results also extend the work of \cite{carlos2}, which examined a three-node bivirus network with a specific tree structure.

To begin, a general convergence result can be established using \Cref{lem:monotone_convergence} if the bivirus system can be shown to have a finite number of equilibria. As we will shortly prove, given generic parameter matrices $D^{i}, B^{i}, i = 1,2$, the bivirus system in \eqref{eq:twovirus} has a finite number of equilibria. We first define what we mean by ```free parameters" and ``generic''. 

\begin{definition}\label{def:free_param}
Free parameters in $B^i,D^i$ are those which are not a priori constrained to assume fixed values for all systems of interest.
\end{definition}

Off-diagonal entries of $D^i$ are not free parameters, being always zero. If $D^i$ is constrained to be the identity matrix (as occurs later), the diagonal entries are also not free parameters. All entries of $B^i$ are free parameters. The free parameters of $D^i$ and $B^i$ take values from the nonnegative real interval, with those of $D^i$ required to be strictly positive while those of $B^i$ are only assumed to be nonnegative.

\begin{definition}
Generic values of the free parameters of $B^i,D^i$ are those not lying on a certain algebraic set of measure zero. 
\end{definition}

At the very least, the existence of such a set would need to be demonstrated to conclude that generic values exist. Often, such a set can be characterized. For the bivirus system, the particular algebraic set is identified as part of the proof of the following result, and is presented in \cref{app:pf_convergence}.

\begin{theorem}\label{thm:convergence}
For generic parameter matrices $D^{i}, B^{i}, i = 1,2$, the bivirus equation set \eqref{eq:twovirus} has a finite number of equilibria. If $D^i=I_n, i=1,2$, then for generic parameter matrices $B^i, i=1,2$ the same conclusion holds. Moreover, for all initial conditions $(x^1(0), x^2(0))\in \Delta$, except possibly for a set of measure zero, the system \eqref{eq:twovirus} will converge to an equilibrium. If the system does not converge to an equilibrium, then it is on a nonattractive limit cycle.
\end{theorem}

This thoroughly answers the two questions posed at the end of \cref{ssec:bivirus}.

% We show that for generic bivirus systems, there are a finite number of equilibria and convergence to an equilibrium with properties satisfying \Cref{lem:interiorequilibria} is the typical outcome (for all initial conditions except possibly a set of measure zero). 

Concerning the first, and as detailed at the end of \cref{ssec:bivirus}, when $\mathcal{R}^i \leq 1$ for some $i = 1,2$, there are at most 3 equilibria; these 3 equilibria remain present when $\mathcal{R}^i > 1$ for $i=1,2$. Sahneh and Scoglio \cite{sahneh2014competitive} showed then that for some parameter choices, there exist \textit{coexistence} equilibria, in which both virus are present. However, \cite{sahneh2014competitive} never proves (or even discusses) if the coexistence equilibrium is unique or indeed whether there are a finite number of them. Simulations suggest that one typically expects a finite number of equilibria, but a special scenario that yields a connected set containing an infinite number of equilibria was identified in~\cite{liu2019bivirus}. 

To keep the exposition clear, the precise algebraic set mentioned above is given in \cref{app:pf_convergence}, and we remark that for any $n$ nodes, the algebraic set is always easily identified. The important point to notice is that the existence of such a set with measure zero (understood in the usual sense for real numbers) establishes the notion of ``generic'' parameter matrices $D^i$ and $B^i$.

% The work of \cite{santos2015bivirus,santos2015bivirus_conference} gives sufficient conditions on $D^i$ and $B^i$ such that one of the boundary equilibria is the global attractor and there are no coexistence equilibria.

Concerning the second, previous works such as \cite{liu2019bivirus,yang2018bivirus,santos2015bivirus,santos2015bivirus_conference,carlos2} establish sufficient conditions on the $D^i, B^i$ parameter matrices to ensure a \textit{specific} equilibrium among the healthy and two boundary equilibria is globally attractive for all initial conditions in $\Delta$. In contrast, \Cref{thm:convergence} establishes that the dynamical system itself is ``almost globally stable'', in the sense that the typical outcome is convergence to an equilibrium, not necessarily the healthy and two boundary ones as illustrated in subsection~\ref{sssec:n2}. This provides a general conclusion on the limiting behaviour of the dynamical system, as opposed to a general condition for a particular equilibrium to be globally attractive. 

Our result does not specify which equilibrium the system converges to and this is deliberate. Identifying conditions on $D^i, B^i$ that yield precise conclusions on the number of equilibria and associated regions of attraction is often difficult~\cite{janson2020networked,santos2015bivirus,santos2015bivirus_conference}. Indeed, there are instances of the bivirus dynamics with no global attractor equilibrium~\cite{carlos2}. In the sequel, we provide an explicit example of such a network, whereby the two boundary equilibria each have a region of attraction with nonzero measure and the boundary between the two regions of attraction forms the stable manifold of an \textit{unstable coexistence equilibrium}. 

The remainder of this paper seeks to explore the properties of the equilibria, and we establish a number of different conclusions i) based on conditions for the $D^i$ and $B^i$ parameter matrices, and ii) exploiting the monotonicity of the bivirus trajectories.

\subsection{Properties of the Equilibria}

Because the Jacobian under transformation in \eqref{eq:jacob_transform} is a Metzler matrix, we can identify different collections of differential equations that give rise to the same equilibria that have the same local stability properties. We start with the following result.

\begin{lemma}\label{prop:Didentity}
Consider two bivirus network systems $\mathcal S$ and $\hat{\mathcal S}$, defined by quadruples $B^1,D^1,B^2,D^2$ and $ \hat B^1=(D^1)^{-1}B^1,\hat D^1= I_n,\hat B^2=(D^2)^{-1} B^2,\hat D^2=I_n$, respectively. Then, the two systems have the same equilibrium sets and the (local) stability properties of each equilibrium are the same. 
% {\color{cyan} If global stability holds for all initial coniditons in $\Delta$ for one system, it also holds for the other system.}
\end{lemma}

This proposition, with proof given in \cref{app:pf_Didentity}, provides an important conclusion: for many theoretical investigations, there will be no loss of generality in assuming $D^1=D^2=I_n$, which can simplify computations and reduce parametrization of the system to just the two matrices $B^1,B^2$ defining the infection rates between nodes. To this end, note that the existence of $(D^i)^{-1}$ ensures that $(D^i)^{-1}B^i$ is irreducible if and only if $B^i$ is irreducible. So the \nameCref{standing_assump} `flows through' if the simplification is undertaken. 
% {\color{red} If the above cyan insertion is made the following blue comes out. }
{\color{black}We further remark that the lemma makes no assertion that global stability (should it be present) in one system implies the same property for the other. Nonetheless, this property is indeed true, as we demonstrate with a simple extension of the lemma in the sequel. }

\begin{remark}
This result allows us to examine scenarios in which the two viruses have different time-scales, by associating $\mathcal S$ and $\hat{\mathcal S}$ in \cref{prop:Didentity} with the quadruples $B^1,D^1,\epsilon B^2,\epsilon D^2$ and $\hat B^1 = (D^1)^{-1}B^1$, $\hat D^1 = I_n$, $\hat B^2 = (\epsilon D^2)^{-1} B^2$, $\hat D^2 = I_n$, respectively, for a positive $\epsilon$ that captures the relative rate of evolution of the two viruses. The location and the \textit{stability} of equilibria remain invariant as $\epsilon$ varies in magnitude. While the invariance of the equilibria locations may seem intuitive, it is certainly not a trivial conclusion that the stability property of each equilibrium is preserved when the relative time-scale of a coupled nonlinear system is changed~\cite{sastry1999nonlinearbook}. On the other hand, varying $\epsilon$ can certainly change the region of attraction of equilibria.
\end{remark}

% The same proof technique used to establish \Cref{prop:Didentity} will also deliver the following result, which roughly states that the location of equilibria and their stability remain invariant if the relative rate of evolution of the two viruses is varied. No actual proof is provided.

% \begin{proposition}
% Consider two bivirus systems $\mathcal S$ and $\overline{\mathcal S}$ defined by quadruples $B^1,D^1,B^2,D^2$ and $B^1,D^1,\epsilon B^2,\epsilon D^2$, respectively, for arbitrary positive $\epsilon$. Then the two systems have the same equilibrium sets and the stability properties of each equilibrium are the same. 
% \end{proposition}

\subsubsection{Nongeneric networks with an infinite number of equilibria}

We will now show that with a generic $B^1$, there are an infinite number of $B^2$, satisfying a specific functional form, which yield a connected set containing an infinite number of equilibria, and this set comprises an interval of a straight line. We term this ``a line of equilibria'' for convenience.  

For convenience, but with no loss of generality as demonstrated by \Cref{prop:Didentity}, take $D^1=D^2=I_n$ and $B^1$ an arbitrary irreducible matrix. We require that $s(-I_n+B^1)>0$, thereby ruling out the possibility of the healthy equilibrium $(\bar x^1,\bar x^2)=({\bf{0}}_n, {\bf{0}}_n)$ being attracting for the bivirus system for any $B^2$ when $(x^1(0), x^2(0)) \in \Delta$. Let ${\bf{1}}_n\gg z\gg{\bf{0}}_n$ with $Z={\rm{diag}}(z)$, satisfy
\begin{equation}\label{eq:zdefinition}
[-I_n+(I_n-Z)B^1]z={\bf 0}_n.
\end{equation}
That is, $z$ is the unique endemic equilibrium of the single virus system \eqref{eq:sis} with $D = I_n$. Now let $C$ be any other nonnegative irreducible matrix for which $z$ is also an eigenvector with unity eigenvalue. Obviously it is straightforward to find such a matrix, and there is an infinite number. The Perron--Frobenius Theorem~\cite{horn2012matrixbook} implies $\rho(C) = 1$ and $z$ is the only positive eigenvector of $C$  (up to a scaling). Define
\begin{equation}\label{eq:B2definition}
B^2=(I_n-Z)^{-1}C.
\end{equation}
\begin{proposition}\label{prop:equilibriumline}
With $D^1=D^2=I_n$, with $B^1$ an arbitrary nonnegative irreducible matrix and with $z$ and $B^2$ defined using \eqref{eq:zdefinition} and \eqref{eq:B2definition}, a set of equilibrium points of the bivirus equilibrium equations \eqref{eq:equil} is given by $\left(\alpha z, (1-\alpha) z\right)$, with $\alpha\in[0,1]$. This equilibrium set is locally exponentially attractive, i.e., for all allowed initial conditions of the bivirus system sufficiently close to the set, the trajectory will approach the set exponentially fast.
\end{proposition}
The proof is given in \cref{app:pf_equibline}.

In \cite{liu2019bivirus}, a line of equilibria was shown to exist if $D^1=k D^2, B^1=kB^2$ for some positive scalar $k$ (and at one point even more specialized conditions). Here, we have relaxed this condition by identifying that a line of equilibria can exist under a broader type of nongenericity; the result in \cite{liu2019bivirus} is covered by \Cref{prop:equilibriumline} after applying \Cref{prop:Didentity}. Whether it is possible to obtain an infinite number of equilibria not lying on a straight line is at the moment an open question. 

It is often appealing to consider special examples of systems where the choice of parameters makes the derivation of analytic results possible. There is a risk however that nongenericity of the examples leads to nongeneric conclusions.  The existence of nongeneric values of matrices defining a bivirus problem  which simplify computations but give rise to an infinite number of equilibrium points is an illustration of the possibility.

In the situation described above, the connected set of equilibria forms an interval of a straight line, and each end of the interval corresponds to each of the single virus equilibria associated with the two viruses being active one at a time. In effect, a bifurcation produces this result. Suppose that $B^2$ above were replaced by $\mu(I_n-Z)^{-1}C$ where $\mu$ is a scalar positive parameter. If $\mu<1$, one can verify that $(\bar x^1, \bar x^2)=(z,{\bf{0_n}})$ is a locally exponentially stable equilibrium, while when $\mu>1$, one can check that it is a saddle point equilibrium, with divergent trajectories. When $\mu=1$, the Jacobian acquires a zero eigenvalue and the bifurcation occurs.   A similar sort of analysis applies if $B^2$ is as above, but $B^1$, after determination of $z$, is then replaced by $\mu B^1$ and one considers then the equilibrium at $(\bar x^1,\bar x^2)=({\bf{0}}_n,z)$ with $\mu$ varying from less than to greater than one. This aspect will now be further explored.    

\subsubsection{Stability of boundary equilibria}
We now explore the stability properties of the boundary equilibria $(\bar x^1, {\bf 0}_n)$ and $({\bf 0}_n,\bar x^2)$ where $\bar x^1 \gg {\bf 0}_n$ and $\bar x^2 \gg {\bf 0}_n$ are separately the unique endemic equilibria for each of the two single virus systems. While $\bar x^1$ and $\bar x^2$ are almost globally stable in the single virus systems, the stability of $(\bar x^1, {\bf 0}_n)$ and $({\bf 0}_n,\bar x^2)$, local or global, in the bivirus case are not guaranteed.

\begin{theorem}\label{thm:boundary_stability}
Consider a generic bivirus system with parameter matrices $D^1 = D^2 = I_n$ and $B^1, B^2$, and suppose that $\rho(B^1) > 1$ and $\rho(B^2) > 1$. Let $\bar x^1\gg {\bf 0}_n$ and $\bar x^2\gg {\bf 0}_n$ denote the equilibria of the two separate single virus systems. Then, 
\begin{enumerate}
    \item The boundary equilibrium $(\bar x^1, {\bf 0}_n)$ is locally exponentially stable if and only if $\rho\big((I_n-\bar X^1)B^2\big) < 1$, and unstable if $\rho\big((I_n-\bar X^1)B^2\big) > 1$.
    \item The boundary equilibrium $({\bf 0}_n, \bar x^2)$  is locally exponentially stable if and only if $\rho\big((I_n-\bar X^2)B^1\big) < 1$, and unstable if $\rho\big((I_n-\bar X^2)B^1\big) > 1$.
\end{enumerate} 
\end{theorem}

The above result, with proof appearing in \cref{app:pf_boundary_stability}, provides necessary and sufficient conditions for local stability of the boundary equilibria, which can be determined by examining two separate single virus systems, with iterative algorithms available for computing $\bar x^1$ and $\bar x^2$, e.g. \cite[Theorem~5]{vanMeighem2009_virus} and \cite[Theorem~4.3]{mei2017epidemics_review}. Some insightful sufficient conditions can be obtained as follows (see \cref{app:pf_boundary}).

\begin{corollary}\label{cor:boundary}
With notation as above, the following statements hold:
\begin{enumerate}
    \item If $B^2 > B^1$, $(\bar x^1, {\bf 0}_n)$ is unstable, and $({\bf 0}_n, \bar x^2)$ is locally stable,  and there is no coexistence equilibrium $(\tilde x^1,\tilde x^2) \in \Delta$.
    \item 
    % Define $\bar b^2 = \max_{i} \sum_{j=1}^n\beta_{ij}^2$ and $\underline{b}^1 = \min_{k} \sum_{j=1}^n \beta_{kj}^1$. If $\bar b^2>\underline{b}^1$,
    If $\underline{b}_2\triangleq \min_{i} \sum_{j=1}^n\beta_{ij}^2 >\bar{b}_1 \triangleq \max_{k} \sum_{j=1}^n \beta_{kj}^1$, then $(\bar x^1, {\bf 0}_n)$ is unstable, and $({\bf 0}_n, \bar x^2)$ is locally stable,  and there is no coexistence equilibrium $(\tilde x^1,\tilde x^2) \in \Delta$.
    \item If $\bar x^2 > \bar x^1$ then $(\bar x^1, {\bf 0}_n)$ is unstable, and $({\bf 0}_n, \bar x^2)$ is locally stable.
\end{enumerate}
\end{corollary}

Item~1 and Item~2 of the corollary give different sufficient conditions on $(\bar x^1, {\bf 0}_n)$ and $({\bf 0}_n, \bar x^2)$ to be unstable and locally stable, respectively. Since there are no coexistence equilibria $(\tilde x^1,\tilde x^2)$, the only other equilibrium is the healthy equilibrium. Later, we present \cref{cor:no_interior_equib}: when there is no coexistence equilibria, one of the two boundary equilibria is in fact stable for all initial conditions in $\Delta$ (and in this instance, it would be $({\bf 0}_n, \bar x^2)$). Notice that Item~1 conditions are entry-wise on the $B^1$ and $B^2$, while Item~2 concerns row sums. Neither subsumes the other, and it is possible to find examples of $B^i$ that satisfy one condition but not the other. Item~2 first appeared in \cite[Corollary~4]{santos2015bivirus_conference} (which refines \cite[Theorem~20]{santos2015bivirus}), and here we give an alternative and short proof. The sufficient condition in Item~3 of the corollary is implied by the condition in Item~1, but the converse is not necessarily true. That is, $B^2 > B^1 \Rightarrow \bar x^2 \gg \bar x^1$, as detailed in \cite[Proposition~4]{liu2019bivirus}. 
% It is interesting to note that Item~1 and Item~2 conditions separately imply that $\mathcal{R}^2 = \rho(B^2) > \mathcal{R}^1 = \rho(B^1)$, as per \cite[Theorem~2.7]{varga2009matrix_book}, providing a link to the reproduction number of the separate single virus cases.

% It is interesting to note that the first sufficient condition in the corollary is straightforward to check, and implies $\mathcal{R}^2 = \rho(B^2) > \mathcal{R}^1 = \rho(B^1)$, as per \cite[Theorem~2.7]{varga2009matrix_book}, {\color{black}allowing a perspective on local stability/instability of boundary equilibria from the reproduction number of the separate single virus cases. Moreover, there are no coexistence equilibria $(\tilde x^1,\tilde x^2)$, and thus only three equilibria for such a bivirus system, viz. the healthy equilibrium, $(\bar x^1, {\bf 0}_n)$, and $({\bf 0}_n, \bar x^2)$, with the former two unstable and the latter locally stable. Later, we present \cref{cor:no_interior_equib}, which then establishes that $({\bf 0}_n, \bar x^2)$ is in fact stable for all initial conditions in $\Delta$. The conditions for global convergence to $({\bf 0}_n, \bar x^2)$ given in \cite[Theorem~20]{santos2015bivirus} and \cite[Corollary~4]{santos2015bivirus_conference} do not subsume the condition $B^2 > B^1$, or vice versa.} The second sufficient condition in the corollary is implied by the first, but the converse is not necessarily true. That is, $B^2 > B^1 \Rightarrow \bar x^2 \gg \bar x^1$, as detailed in \cite[Proposition~4]{liu2019bivirus}. 

All three configurations of boundary equilibria stability properties can occur: i) both boundary equilibria are locally exponentially stable, ii) both boundary equilibria are unstable, and iii) one boundary equilibria is locally stable and the other unstable. In \cref{sssec:n2} below, we give examples of each configuration for an $n=2$ network.

Knowing the stability configuration of the two boundary equilibria, one can exploit the properties of monotone systems to obtain simple counting results that lower bound the number of interior equilibria, see \cite[Theorem~2.8]{smith1988monotone_survey}. When there are no interior equilibria, a global stability result for one of the boundary equilibria follows. We explore this in further detail in the next subsection.

\subsection{Properties of trajectories and implications for equilibria}\label{ssec:properties}

% \subsection{Trajectory properties arising from the monotone system property}

In this subsection, we argue that for a given system, all trajectories with $x^i(0) \gg {\bf 0}_n, i = 1,2$ are bounded by two special trajectories, starting from special corners of the allowed set of initial conditions. Other trajectory bounding results have appeared in \cite{santos2015bivirus,santos2015bivirus_conference}, but by using two special trajectories, we are able to go further by demonstrating that from this result flows several important conclusions on equilibria properties. For instance, the results provide a simulation tool involving construction of two trajectories which on occasions can be expected to exclude the possibility of any other stable equilibria than those associated with the two trajectories, and on other occasions, to flag the presence of an unstable equilibrium. We draw on the conclusion, recorded after \Cref{lem:jacobianirreducible}, that the bivirus system is type $K_m$ monotone with $m=({\bf 0}_n^\top,{\bf 1}_n^\top)$.

\begin{theorem}\label{thm:trajectoryordering}
Consider the equation set \eqref{eq:twovirus}, and in particular consider the trajectories $x_A(t),x_B(t)$ defined for arbitrarily small but positive $\eta$ by the initial conditions $x_A^1(0)=\frac{1}{2}\eta{\bf{1}}, x_A^2(0)=(1-\eta){\bf{1}}$ and $x_B^1(0)=(1-\eta){\bf{1}}, x_B^2(0)=\frac{1}{2}\eta{\bf{1}}$. Suppose $x_C(t)$ is a trajectory beginning at any initial condition satisfying $x^1_B(0)> x^1_C(0)> x^1_A(0)$, and $x^2_A(0)> x^2_C(0) > x^2_B(0)$, and $x_C^1(0)+x^2_C(0)< {\bf{1}}_n$.
% \begin{eqnarray}
% x^1_B(0)> x^1_C(0)> x^1_A(0)\\\nonumber
% x^2_A(0)> x^2_C(0) > x^2_B(0)\\\nonumber
% x_C^1(0)+x^2_C(0)< {\bf{1}}_n
% \end{eqnarray}
Then the following inequalities hold for all $t>0$
\begin{align}\label{eq:3trajectory}
  x_B^1(t)\gg x_C^1(t)\gg x_A^1(t)\,,\quad 
  x_A^2(t)\gg x_C^2(t)\gg x_B^2(t)\,,\quad 
  x_C^1(t)+x_C^2(t)\ll {\bf {1}}_n.
\end{align}
% \begin{align}\label{eq:3trajectory}
%   x_B^1(t)\gg x_C^1(t)&\gg x_A^1(t)\\\nonumber
%   x_A^2(t)\gg x_C^2(t)&\gg x_B^2(t)\\\nonumber
%   x_C^1(t)+x_C^2(t)&\ll {\bf {1}}_n
% \end{align}
\end{theorem}

\Cref{app:pf_trajorder} provides a very short proof based on known properties from the monotone systems literature, while \cite{santos2015bivirus} obtains a similar result using computations tailored to the bivirus dynamics. 

% {\color{red} If the cyan extemsion to Lemma 3.7 is accepted, the following two blue paragraphs including the new lemma come out.}
{\color{black}The theorem also allows us to directly extend \Cref{prop:Didentity} to deal with bivirus systems with a globally stable equilibrium, as stated in the following result (with proof in \Cref{app:pf_global_stab}). 

\begin{lemma}\label{lem:global_stab}
Assume the same hypotheses as \Cref{prop:Didentity}, with the two systems $\mathcal{S}$ and $\hat{\mathcal{S}}$. If one system has an equilibrium which is globally stable for all initial conditions in $\Delta$, then the global stability property holds for the second system.
\end{lemma}
}
We further expand the theorem via a second corollary dealing with the limit points of the trajectories $x_A(t), x_B(t), x_C(t)$. Before stating it, we make an important observation, based on \Cref{lem:monotone_convergence} and the discussions below it. Specifically, we know that
% It is known from \cite{smith1988monotone_survey} that if the set of equilibria is countable, a property which is fulfilled by generic $D^i,B^i$ as detailed in \Cref{prop:finite_solution_set}, then for every initial condition in the interior of $\Xi_n\times \Xi_n$ other than those on a set of measure zero, the trajectory will converge to a  stable equilibrium point. From an initial condition in the exceptional set of measure zero, convergence must occur to a saddle point, or the initial condition lies on a nonattractive limit cycle. 
if convergence to a stable equilibrium point does not occur for one of the initial conditions used in \Cref{thm:trajectoryordering}, then convergence will occur to a stable equilibrium for a perturbation of that initial condition within a ball of arbitrarily small radius. We abbreviate this idea below with the words `perturbed if necessary'. 
\begin{corollary}\label{cor:equilibriumordering}
With the same hypothesis as \Cref{thm:trajectoryordering},  with values for the matrices $D^i,B^i$  assuring the number of equilibria is finite,  and with initial conditions perturbed if necessary,  assume that the trajectories $x_A(t),x_B(t)$ approach limits $\bar x_A,\bar x_B$ respectively. If $\bar x_A\neq \bar x_B$, let $\mathcal W = \{x : x_A \leq_{K_m} x \leq_{K_m} x_B\}$, with $m = ({\bf 0}_n^\top,{\bf 1}_n^\top)$ as above, denote the closed hyperrectangle\footnote{In some cases, one or more coordinates of $\bar x_A$ and $\bar x_B$ may assume the same value, (though if $\bar x_A\neq \bar x_B$ neither $\bar x_{A}^{1}=\bar x_{B}^{1}$ nor $\bar x_{A}^{2}=\bar x_{B}^{2}$ is possible by \Cref{lem:interiorequilibria}), and  $\mathcal{W}$ is then degenerate.} whose edges are axis-parallel and which has a principal diagonal joining the points $\bar x_A$ and $\bar x_B$. Then
\begin{enumerate}
    \item 
If $x_C(0)$ lies outside $\mathcal W$,  the associated trajectory converges to a limit $\bar x_C$ within $\mathcal W$, obeying the constraints  implied by Theorem \ref{thm:trajectoryordering}.
\item 
If $x_C(0)$ lies within  $\mathcal W$, the associated trajectory either converges to a limit $\bar x_C$ again obeying the constraints implied by Theorem \ref{thm:trajectoryordering}  or lies on a nonattractive limit cycle, which requires that $\mathcal W$ has dimension at least 2.
\item
If $\bar x_A$ and $\bar x_B$ coincide, then $\bar x_C=\bar x_A=\bar x_B$.
\end{enumerate}
\end{corollary}

\begin{figure}
    \centering
    \def\svgwidth{0.6\linewidth}
    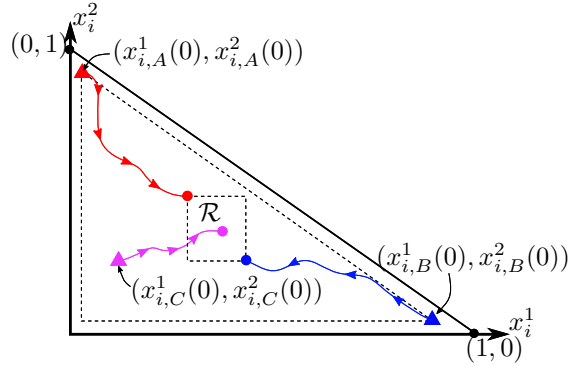
    \caption{Illustration of Corollary~\ref{cor:equilibriumordering}. The triangles indicate initial conditions, with the system trajectory shown converging to the circles, which indicate the equilibria. The dotted rectangle shows the intersection of the hyperrectangle $\mathcal{W}$ with the plane defined by the $i$th coordinates of $x^1$ and $x^2$.} 
    \label{fig:bivirus_trajectories}
\end{figure}

The above corollary is illustrated in \Cref{fig:bivirus_trajectories}, and proved in \cref{app:pf_equiborder}. The corollary effectively offers  a simulation-type test for establishing whether there is a single equilibrium for which $\bar x^i\gg {\bf 0}_{n}, i=1,2$: one simply computes two trajectories, $x_A(t)$ and $x_B(t)$, from initial conditions as close as possible to opposite corners of the region of interest, and checks that they approach a common limit. What happens if the two trajectories do not approach a common limit? The answer, again appealing to monotone systems theory, is that an unstable equilibrium must lie in $\mathcal{W}$. 
% If we are able to identify parameters yielding two such trajectories, we would prove that it is possible to have an unstable coexistence equilibrium; currently, this is an open problem. 
% In fact,  we do not know (at this stage)  whether it is ever possible to have an unstable (including saddle point) equilibrium $\bar x_C$ satisfying ${\bf{1}}_{2n}\gg \bar x_C\gg {\bf 0}_{2n}$, though certainly monotone system theory does \textit{not} exclude such a possibility. 
The following is an amalgam of Theorem 2.8 and Proposition 2.9 in \cite{smith1988monotone_survey}, and thus no proof is presented.
% , a condition for existence can be found, and there are  several properties of the unstable trajectories associated with it. 

\begin{corollary}
For bivirus system \eqref{eq:twovirus} and with generic parameter matrices, suppose there are two attractive equilibria $\bar x_A$ and $\bar x_B$ with $\bar x_A\gg_{K_{m}}\bar x_B$. Then  there exists an unstable equilibrium $\bar x_C$ obeying
\begin{equation}\label{eq:threeequilibria}
\bar x_A\gg_{K_{m}}\bar x_C\gg_{K_{m}}\bar x_B,
\end{equation}
and conversely, if an unstable equilibrium $\bar x_C$ exists with $\bar x_C\gg_{K_{m}} {\bf 0}_{2n}$, there exist two attractive equilibria $\bar x_A$ and $\bar x_B$ satisfying \eqref{eq:threeequilibria}. 
Further, if there is no other stable equilibrium $\bar x_D$ satisfying $\bar x_A\gg \bar x_D\gg \bar x_B$, then there are two unstable trajectories emanating from the equilibrium $\bar x_C$ that are tangent to the eigenvector associated with the unstable eigenvalue of the associated Jacobian matrix $J(\bar x_C)$, satisfying the following properties. The first trajectory satisfies $\dot x^1\gg {\bf 0}_{n}$ and $-\dot x^2\gg {\bf 0}_{n}$ along its entirety and its limiting point is $\bar x_B$. The second trajectory satisfies $-\dot x^1\gg {\bf 0}_{n}$ and $\dot x^2\gg {\bf 0}_{n}$ along its entirety and its limiting point is $\bar x_A$.
% \begin{equation}
% \dot x^1\gg {\bf 0}_{n},-\dot x^2\gg {\bf 0}_{n}
% \end{equation}
% along its entirety and its limiting point is $\bar x_B$. The second trajectory satisfies
% \begin{equation}
% -\dot x^1\gg {\bf 0}_{n}, \dot x^2\gg {\bf 0}_{n}
% \end{equation}
% along its entirety and its limiting point is $\bar x_A$.
\end{corollary}

Finally, if the particular system has no coexistence equilibria, then, as recorded in the next corollary the properties established above imply that one of the boundary equilibria is globally stable, being attractive for all initial conditions in $\Delta$. Sufficient conditions on $D^i, B^i$ for there to be no coexistence equilibria include \cref{cor:boundary} and \cite[Theorem~6]{janson2020networked}. The proof is in \cref{app:pf_no_interior_equib}.

\begin{corollary}\label{cor:no_interior_equib}
For the bivirus system in \eqref{eq:twovirus} with generic parameter matrices, suppose there are no coexistence equilibria. Then precisely one of the boundary equilibria is an attractive equilibrium, and the set $\Delta$ is in its region of attraction.
\end{corollary}

\subsection{Equilibria for a low order system}\label{sssec:n2}
We first explain how, for generic values of system parameters and with $n=2$, one can compute using nothing more than a single quadratic equation solution followed by linear equation solutions, any equilibrium for which $\bar x^1\gg{\bf 0}_{n}, \bar x^2\gg {\bf 0}_{n}$. Then we illustrate a series of possible outcomes which provides a comprehensive account of the possible limiting behaviours for the bivirus system. 
% There is an important consequence. Because of the occurrence of this quadratic equation in the solution procedure for $n=2$, there can be at most two equilibria lying in the interior of the region of interest. (Note that the quadratic equation in principle can have real solutions outside the region of interest, or imaginary solutions) The previous results of \cref{ssec:properties} show that 
%  if there are two distinct stable equilibria (either one or both being permitted to be on a boundary) , there would have to be  a third saddle point or unstable equilibrium `between' them. There cannot be simply two equilibria with all entries positive.  However, because of the occurrence of a quadratic equation in the solution procedure when $n=2$, there can be at most two equilibria (with no account taken of where they lie) found by the procedure, and so there can actually be no more than one equilibrium with $\bar x^1\gg{\bf 0}_{2}, \bar x^2\gg {\bf 0}_{2}$ and it must be stable.  Note that by \Cref{lem:interiorequilibria}, for the systems of interest at any equilibrium for which $\bar x^i\geq {\bf 0}_{2}$, either $\bar x^{i}={\bf 0}_{2}$ or $\bar x^{i}\gg {\bf 0}_{2}$.

The solution procedure is straightforward. Without loss of generality, assume $D^1=D^2=I_2$. Setting $\alpha=(\bar x^1_2/\bar x^1_1)$ and $\gamma=(\bar x^2_2/\bar x^2_1)$, the equilibrium equations in \eqref{eq:equil} yields after some manipulation
\begin{align*}
\beta_{11}^1+\beta^1_{12}\alpha=\beta^2_{11}+\beta^2_{12}\gamma\,,\quad \text{and} \qquad \beta^1_{21}\alpha^{-1}+\beta^1_{22}=\beta^2_{21}\gamma^{-1}+b^2_{22}
\end{align*}
% Then the equilibrium equations can be recorded as follows:
% \begin{align}
%     -\bar x^1_1+(1-\bar x^1_1-\bar x^2_1)(b^1_{11}\bar x^1_1+b^1_{12}\bar x^1_2)&=0\\\notag
%     -\bar x^1_2+(1-\bar x^1_2-\bar x^2_2)(b^1_{21}\bar x^1_1+b^1_{22}\bar x^1_2)&=0\\\notag
%     -\bar x^2_1+(1-\bar x^1_1-\bar x^2_1)(b^2_{11}\bar x^2_1+b^2_{12}\bar x^2_2)&=0\\\notag
%     -\bar x^2_2+(1-\bar x^1_2-\bar x^2_2)(b^2_{21}\bar x^2_1+b^2_{22}\bar x^2_2)&=0
% \end{align}
% Setting $\alpha=(\bar x^1_2/\bar x^1_1)$ and $\beta=(\bar x^2_2/\bar x^2_1)$ yields 
% \begin{align*}
% b_{11}^1+b^1_{12}\alpha&=b^2_{11}+b^2_{12}\beta\\ b^1_{21}\alpha^{-1}+b^1_{22}&=b^2_{21}\beta^{-1}+b^2_{22}
% \end{align*}
One can eliminate $\gamma$ and obtain a quadratic equation for $\alpha$. For each solution of the quadratic equation, it is then possible using only linear equations to obtain values for the $\bar x^i_j$. Note that these have no constraints on the signs of their entries. The fact that a quadratic equation underpins the algortihm means that there can be at most two interior equilibria. Exploiting properties of monotone systems, an argument centred around \cite[Proposition~2.9]{smith1988monotone_survey} will lead to the conclusion that \textit{if} there are two interior equilibria, one must be stable and the other unstable. In the examples below, we report either no interior equilibria, or the interior equilibrium is unique.

% We have not seen an example with more than one. When the boundary equilibria are both stable, a single interior equilibrium must be unstable, and when the boundary equilibria are both unstable, a single interior equilibrium must be stable, as will be seen in the discussion of {\color{red} cross reference to section 3.3. of the paper. Ben may need to pick up this point in a minor revision of section 3.3}
% The argument above shows that there can be only one which has $\bar x^1\gg{\bf 0}_{2}, \bar x^2\gg {\bf 0}_{2}$.

% We now show using a specific example the diverse range of outcomes possible in the bivirus dynamics, connections with the results presented above, and the strength of the above solution procedure. 
With $D^1 = D^2 = I_2$, we fix $$B^1=\begin{bmatrix} 1.6&1\\1&1.6\end{bmatrix}$$
and vary the $B^2$ matrix according to \cref{tab:n2_example}. The equilibria $(x^1, x^2)$, with $x^1 = [x^1_1, x^1_2]^\top$ and $x^2 = [x^2_1, x^2_2]^\top$, for each case as computed using the above solution procedure are reported, excluding the healthy equilibrium $({\bf 0}_2, {\bf 0}_2)$.

\textit{Case 1:} A line of coexistence equilibria (which is locally exponentially attractive by \cref{prop:equilibriumline}) joins the boundary equilibrium $(0.615\cdot{\bf{1}}_2,{\bf{0}}_2)$ to the other boundary equilibrium $({\bf{0}}_2, 0.615\cdot{\bf{1}}_2)$. There are no other coexistence equilibria other than those in the line, and extensive simulations suggest the line of equilibria is in fact globally attractive for $\Delta$.
    
\begin{figure*}
\centering
\begin{minipage}{0.49\linewidth}
\centering
\subfloat[Node 1]{\includegraphics[height=\linewidth,angle=-90]{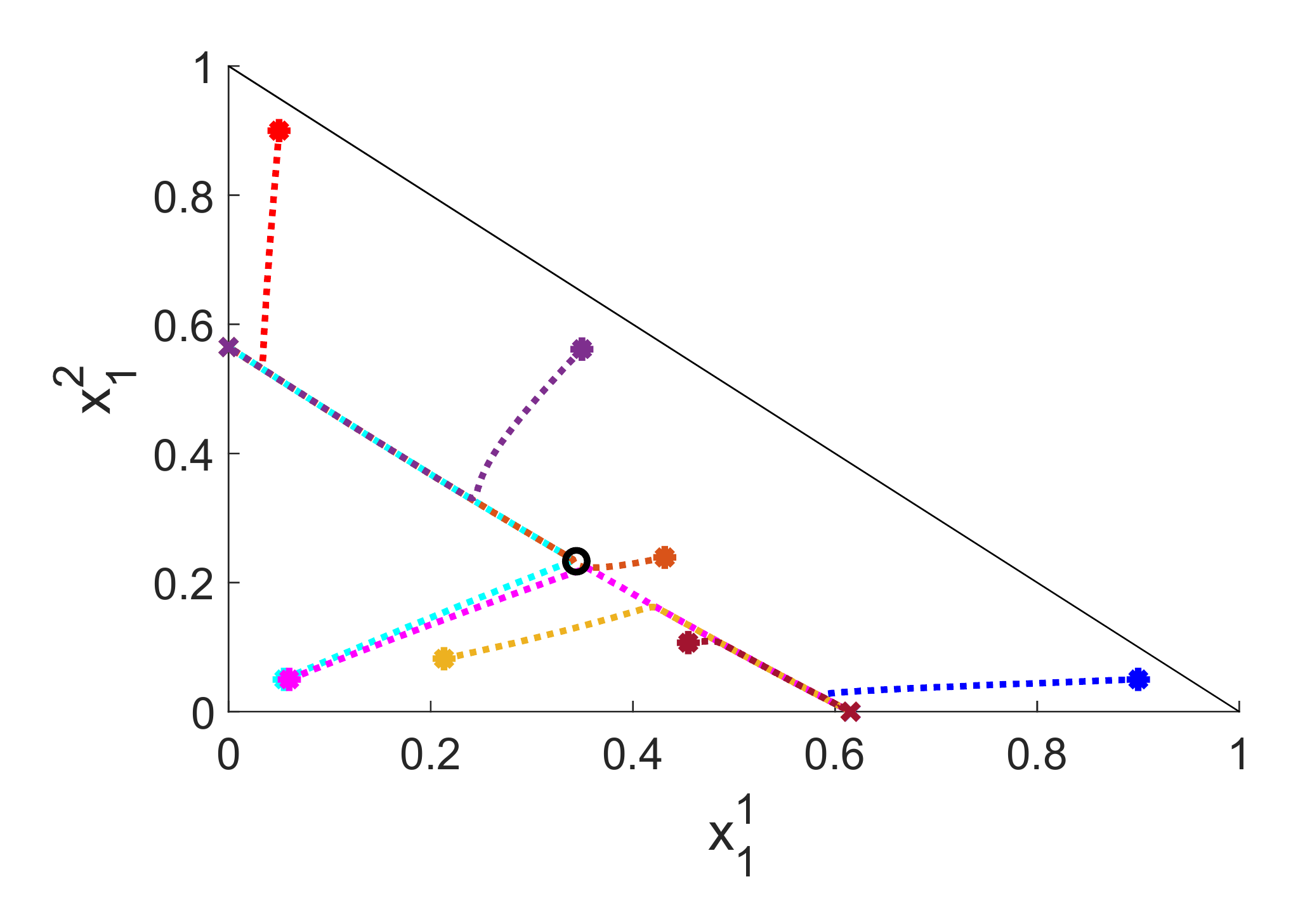}\label{fig:n2_node1}}
\end{minipage}
\hfill
\begin{minipage}{0.49\linewidth}
\centering
\subfloat[Node 2]{\includegraphics[height=\linewidth,angle=-90]{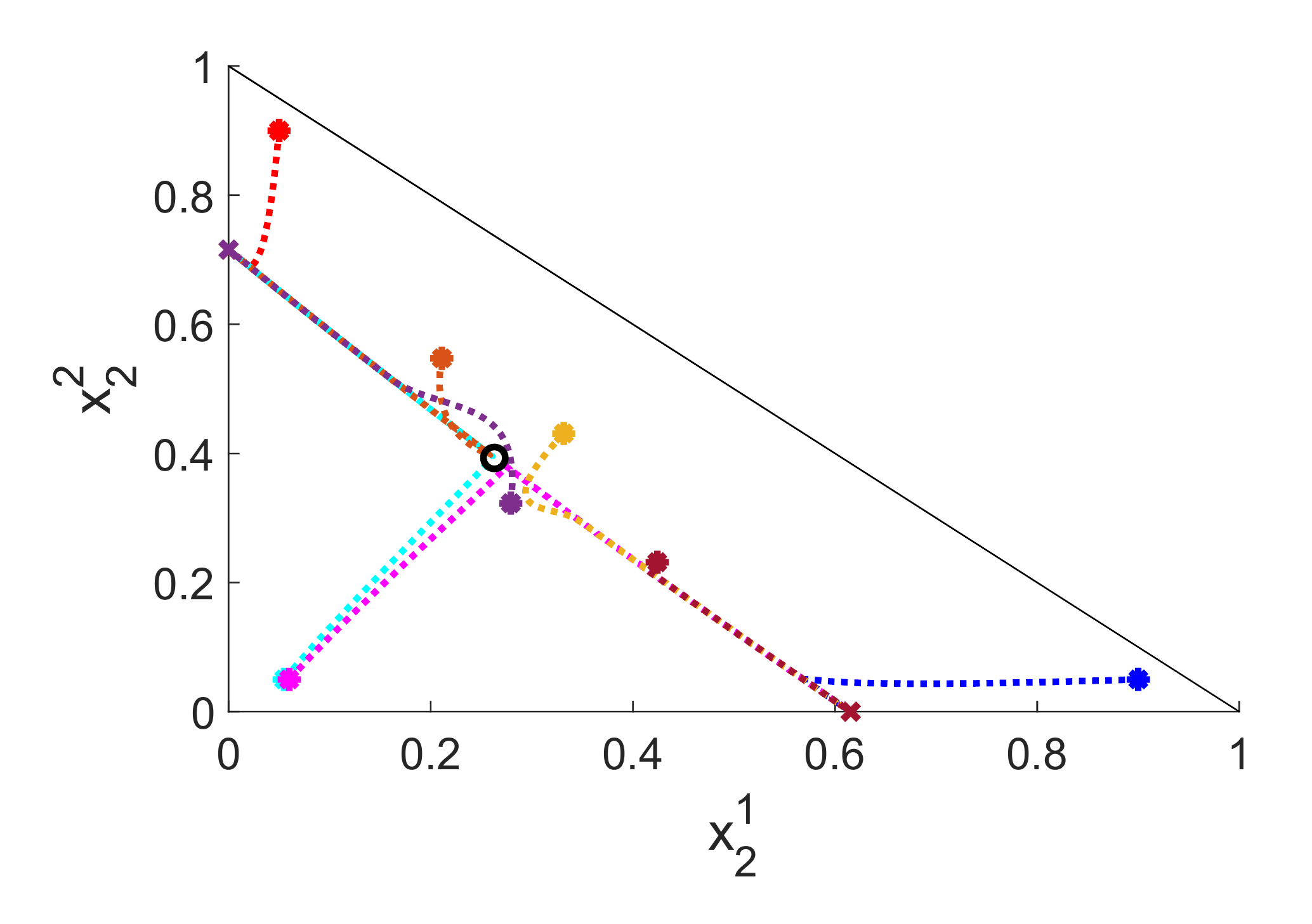}\label{fig:n2_node2}}
\end{minipage}
\caption{Simulation of the bivirus system for \textit{Case 2} identified in \cref{sssec:n2}. Different colored lines represent the trajectories for eight different initial conditions for (a) Node 1 states $x^1_1(t)$ and $x^2_1(t)$ and (b) Node 2 states $x^1_2(t)$ and $x^2_2(t)$. The solid ball denotes the initial condition, and the cross denotes the equilibrium reached as $t\to\infty$. The black circle identifies the unstable coexistence equilibrium $(\tilde x^1, \tilde x^2)$.
% , where $\tilde x^1=[0.344\;0.263]^{\top}, \tilde x^2=[0.233\;0.393]^{\top}$.
}
\label{fig:n2_example}
\end{figure*}

\textit{Case 2:} This example yields two stable boundary equilibria at $(\bar x_1,{\bf{0}}_2)$ and $({\bf{0}}_2, \bar x_2)$ with $\bar x_1=[0.615\;0.615]^{\top}$ and $\bar x^2=[0.565\;0.715]^{\top}$. There is a unique coexistence equilibrium $(\tilde x^1, \tilde x^2)$, and it is unstable: $\tilde x^1=[0.344\;0.263]^{\top}, \tilde x^2=[0.233\;0.393]^{\top}$. In fact, $(\tilde x^1, \tilde x^2)$ is hyperbolic, and by combining properties of attractive manifolds for unstable hyperbolic equilibria of monotone systems, \cite[Theorem~2.10]{smith1988monotone_survey}, and standard results on regions of attraction for stable equilibria~\cite{chiang2015stability}, it follows that the regions of attraction for $(\bar x_1,{\bf{0}}_2)$ and $({\bf{0}}_2, \bar x_2)$ encompass all of $\Delta$ except for a set of measure zero, being the attractive manifold of $(\tilde x^1, \tilde x^2)$ on which there are no two distinct points $y,z$ with $y <_{K_m} z$. In fact, this manifold is part of the boundary of the regions of attraction for $(\bar x_1,{\bf{0}}_2)$ and $({\bf{0}}_2, \bar x_2)$. Fig.~\ref{fig:n2_example} demonstrates this via sample trajectories.
% {\color{red} Ben to fix earlier part of main manuscript in last para of subsection 3.2.2 saying we do not know if this possibility is achievable, and maybe highlight somewhere, here begin just one possibility, the novelty in this conclusion}

\textit{Case 3:} There are two unstable boundary equilibria $(\bar x_1,{\bf{0}}_2)$ and $({\bf{0}}_2, \bar x_2)$ with $\bar x_1=[0.615\;0.615]^{\top}$ and $\bar x^2=[0.665\;0.515]^{\top}$. There is a unique coexistence equilibrium $(\tilde x^1, \tilde x^2)$, and it is locally stable: $\tilde x^1=[0.462\;0.512]^{\top}$ and $\tilde x^2=[0.168\;0.089]^{\top}$. \cite[Proposition~4]{ye2021_PH_TAC} establishes global convergence to $(\tilde x^1, \tilde x^2)$ for all $(x^1(0),x^2(0))\in\Delta$.

\textit{Case 4:} The boundary equilibria $(\bar x_1,{\bf{0}}_2)$ and $({\bf{0}}_2, \bar x_2)$ are unstable and locally stable, respectively: $\bar x_1=[0.615\;0.615]^{\top}$ and $\bar x^2=[0.665\;0.655]^{\top}$. There are no coexistence equilibrium, and \cref{cor:no_interior_equib} thus establishes that $({\bf{0}}_2, \bar x_2)$ is attractive for all initial conditions in $\Delta$.

\begin{remark}
{\color{black}Case 2 has a particularly notable outcome, because in the context of a ``survival-of-the-fittest'' battle between the two viruses, different initial conditions in $\Delta$ can result in either virus surviving. A similar outcome was presented in \cite{carlos2}, but for a three-node network with specific tree structure. Our recent work has proposed a method to systematically construct bivirus networks, with an arbitrary number of nodes and arbitrary topology, where either virus can win the battle~\cite{ye2021_bivirus_outcomes}. Such outcomes are distinct from existing results such as \cite{liu2016bivirus,sahneh2014competitive,santos2015bivirus,santos2015bivirus_conference}, which identify parameter regimes for a specific virus to win the battle independent of initial conditions.}
% Finally, this further underlines the strength of \cref{thm:convergence}, which guarantees convergence to an equilibrium is the typical outcome, even when there is no global attractor equilibrium.
\end{remark}

\begin{table}
    \centering
    \caption{The $B^2$ matrices corresponding to different cases in the $n = 2$ example.}
    \label{tab:n2_example}
\scalebox{0.9}{   \begin{tabular}{c | c | c | c }
Case 1 & Case 2 & Case 3 & Case 4 \\
$B^2 = \begin{bmatrix}2.1&0.5\\1.5&1.1\end{bmatrix}$ & $B^2=\begin{bmatrix} 2.1&0.156\\3.0659&1.1\end{bmatrix}$ & $B^2=\begin{bmatrix}2.1&1.143\\0.745&1.1\end{bmatrix}$ & $B^2=\begin{bmatrix}2.1&0.885\\1.885&1.1\end{bmatrix}$
    \end{tabular}}
\end{table}

\section{Conclusions}\label{sec:conclusions}
In this paper, we have analyzed a deterministic networked bivirus model. We have proved the bivirus system with generic parameters has a finite number of equilibria. Using monotone systems theory, and assuming generic parameters, we established convergence to an equilibrium for (almost) all initial conditions. The properties of equilibria have been further explored, using a mixture of monotone systems theory, matrix theory, and algebraic analysis. Future work will focus on the implications of stable/unstable boundary equilibria on the number and stability of interior equilibria, and control strategies that leverage one virus to eliminate the other, perhaps by optimal design of parameters to ensure a specific boundary equilibrium is globally stable.

\appendix

\section{Proofs of main results} 

\subsection{Proof of \Cref*{lem:interiorequilibria}}\label{app:pf_interiorequib}

As we know, there are equilibria $({\bf 0}_n, {\bf 0}_n), (\bar x^1, {\bf 0}_n)$ and $({\bf 0}_n,\bar x^2)$ for vectors $\bar x^1,\bar x^2$ which are separately equilibria of single virus systems. Because $\bar x^1 \gg {\bf 0}_n$ is the unique endemic equilibrium associated with the dynamics \eqref{eq:sis} if only virus 1 is considered, it is clear that there cannot exist an equilibrium $(\hat x^1, {\bf 0}_n)$ with ${\bf 0}_n \leq \hat x^1 \neq \bar x^1$. By the same reasoning, there cannot exist an equilibrium $({\bf 0}_n, \hat x^2)$ with ${\bf 0}_n \leq \hat x^2 \neq \bar x^2$.

If and only if $\mathcal{R}^i > 1$, or $s(-D^i+B^i)>0$ for $i=1,2$, it becomes possible for further equilibria to exist with $\bar x^i> {\bf 0}_n$. To complete the first part of the proof, we show that there can be no equilibrium $(\bar x^1,\bar x^2)$ in which $\bar x^i\neq {\bf 0}_n$, but $\bar x^i\gg {\bf 0}_n$ fails. 

Observe first that there exists no $j$ for which $1-x^1_j-x^2_j=0$. For if there were such a $j$, the $j$-th row of the two equilibrium equations \eqref{eq:equil} would imply $x^1_j=x^2_j=0$, a contradiction. Hence  $\bar x^1+\bar x^2\ll{\bf 1}_n$ and $I_n-\bar X^1-\bar X^2$ is a nonsingular positive diagonal matrix, so that $(I_n-\bar X^1-\bar X^2)B^i$ and indeed $(D^i)^{-1}(I_n-\bar X^1-\bar X^2)B^i$ for $i=1,2$  are then both irreducible and nonnegative. This means that for any $y> {\bf 0}_n$, there is at least one positive entry of the vector $(D^i)^{-1}(I_n-\bar X^1-\bar X^2)B^iy$ in those rows where $y$ has a zero entry (see \cref{ssec:linear_algebra}). Since however $\bar x^i=(D^i)^{-1}(I_n-\bar X^1-\bar X^2)B^i\bar x^i$, identifying $y$ with $\bar x^i$  would yield a contradiction, unless $\bar x^i\gg {\bf 0}_n$. 

Lastly, observe that $$\bar x^1=(D^1)^{-1}(I_n-\bar X^1-\bar X^2)B^1\bar x^1=(D^1)^{-1}(I_n-\bar X^1-\tilde X^2)B^1\bar x^1,$$ and because $\bar x^1\gg {\bf 0}_n$ and $B^1\bar x^1\gg {\bf 0}_n$ it follows that $\bar X^2=\tilde X^2$ or $\bar x^2=\tilde x^2$ as required.

\subsection{Proof of \Cref*{lem:strictinvariant}}\label{app:pf_strictinvariant}

Set $z={\bf 1}_n-x^1-x^2$ and $\tilde B^i(x^i)={\rm{diag}}(B^ix^i)$. Observe then that $(I_n-X^1(t)-X^2(t))B^ix^i(t)=\tilde B^i(x^i(t))z$.

It is straightforward to derive that 
\begin{align}
 \dot x^i&=-D^ix^i(t)+\tilde B^i(x^i(t))z(t)\,,\quad i = 1, 2\\\notag
%  \dot x^2&=-D^2x^2(t)+\tilde B^2(x^2(t))z(t)\\\notag  
 \dot z&=D^1x^1(t)+D^2x^2(t)-[\tilde B^1(x^1(t))+\tilde B^2(x^2(t))]z(t)
\end{align}
Consider the second equation.  If $z_i(t)=0$ for some $i$, $1-x^1_i(t)-x^2_i(t)=0$, so that one at least of $x^1_i(t)$ and $x^2_i(t)$ is nonzero, i.e. $\dot z_i$ is positive. In fact, there exists $\epsilon$, such that if $z_i(t)\in[0,\epsilon]$, then $\dot z_i(t)>0$; this implies that $z_i(t)$ can never be zero for $t>0$, i.e. $z(t)\gg {\bf 0}_n\;\forall t>0$. Obviously this also establishes that $x^i(t)\ll {\bf 1}_n$ for $t>0$. 

Next, suppose that at some time $T$, $k$ entries of $x^1(t)$ in positions, $i_1,i_2,\dots i_k$ say, are zero, with $k < n$. Then since $B^1$ is irreducible, $B^1x^1$ must have at least one entry in one of the  positions $i_1,i_2,\dots i_k$ which is nonzero, as detailed in \cref{ssec:linear_algebra}. Suppose it is position $i_1$. Then from the above equations,  $\dot x^1_{i_1}(t)>0$.  It follows that for $t>T$ with $t-T$ sufficiently small, there are fewer than $k$ entries of $x^i$ which are zero. An extension of this argument in fact shows that for all $t>0$, all entries of $x^i$ must be nonzero, i.e. $x^i(t)\gg {\bf 0}_n$ for $t>0$.  

\subsection{Proof of \Cref*{lem:jacobianirreducible}}\label{app:pf_jacobirred}

 Under the hypothesis of the lemma, \Cref{lem:strictinvariant} guarantees that $I_n-X^1-X^2$ is nonsingular, and so $(I_n-X^1-X^2)B^i$ is irreducible. Moreover, $\tilde B^i, i=1,2$ are nonsingular positive diagonal matrices since $x^i\gg {\bf 0}_n$ for $i=1,2$. Evidently, it is enough then to show that the nonnegative matrix 
\begin{equation}
    H:=\begin{bmatrix}
    B^1& \tilde B^1\\
    \tilde B^2&B^2
    \end{bmatrix}
\end{equation}
is irreducible, which is equivalent to the graph $G = (\mathcal{V}, \mathcal{E}, H)$ being strongly connected. The graph $G$ is comprised of two strongly connected subgraphs, call them $\mathcal{G}^1$ and $\mathcal{G}^2$, associated with irreducible matrices $B^1$ and $B^2$, respectively. The matrix $\tilde B^1$ captures edges from nodes in $\mathcal{G}^2$ to nodes in $\mathcal{G}^1$, while $\tilde B^2$ captures edges from nodes in $\mathcal{G}^1$ to nodes in $\mathcal{G}^2$. Since $\tilde B^i$ are positive diagonal for $i=1,2$, there exists a path from any node $k$ in $\mathcal{G}^1$ to any node $j$ in $\mathcal{G}^2$, and vice versa. It follows that $G$ is strongly connected, and hence $H$ is irreducible. The proof for the equilibrium value $PJ(\bar x^1,\bar x^2)P$ is effectively identical. 

\subsection{Proof of \Cref*{prop:Didentity}}\label{app:pf_Didentity}
When the two equations of \eqref{eq:equil} corresponding to $\mathcal S$ are multiplied on the left by $(D^1)^{-1}$ and $(D^2)^{-1}$, respectively, equilibrium equations for $\hat{\mathcal S}$ result (and the converse is trivial). The proves that the equilibrium sets are identical. 
The equivalence of equilibrium stability properties requires slightly more work. Begin with system $\mathcal S$.  Because $\tilde J_{\mathcal S}=:PJ_{\mathcal S}(\bar x^1,\bar x^2)P$ is Metzler, local exponential stability of the equilibrium $(\bar x^1,\bar x^2)$ implies that $-\tilde J_{\mathcal S}$ is an $M$-matrix. As the $M$-matrix property is preserved under multiplication by a positive diagonal matrix, it follows on multiplying $\tilde J_{\mathcal{S}}$ by ${\rm{diag}}\left((D^1)^{-1},(D^2)^{-1}\right)$ that the result is another $M$-matrix. But this multiplication yields $\tilde J_{\hat{\mathcal S}}=-PJ_{\hat{\mathcal S}}P$. Thus the stability property for the equilibrium is the same. An identical argument works if $-\tilde J_{\mathcal S}$ is a singular $M$-matrix and by negation if it is not an $M$-matrix. 

% {\color{cyan} Without loss of generality, suppose that $\mathcal{S}$ has a globally stable equilibrium at $\bar x$. To obtain a contradiction, suppose there is an initial state $x(0)\in\Delta$ for $\hat{\mathcal{S}}$ for which the trajectory $x(t)$ does not converge to $\bar x$. Such a point $x(0)$ could not be a coexistence equilibrium point, as it would also be a coexistence equilibrium point in $\mathcal S$, contradicting the global stability assumption. Then, $x(0)$ can only be a point on the stable manifold of an unstable equilibrium point in $\hat{\mathcal S}$. However, the only unstable equilibria of $\mathcal S$, given global stability of $\bar x$ in $\Delta$, are the healthy equilibrium and (possibly) boundary equilibria with either $x^i = {\bf 0}_n$. The same must be true of $\hat{\mathcal{S}}$. But we know that all such trajectories forming the stable manifold are confined to the boundary (where $x^i = {\bf 0}_n$ for one of $i=1,2$), moving from the healthy equilibrium to a boundary equilibrium. This contradicts the fact that $x(0)\in\Delta$. }

\subsection{Proof of \Cref*{prop:equilibriumline}}\label{app:pf_equibline}
Observe that 
\[
[-I_n+(I_n-Z)B^2]z=[-I_n+C]z={\bf 0}_n
\]
Also, setting $\bar x^1_{\alpha}=\alpha z, \bar x^2_{\alpha}=(1-\alpha)z$ yields $Z=\bar X^1_{\alpha}+\bar X^2_{\alpha}$, with $Z$ independent of $\alpha$. Hence
$
[-I_n+(I_n-\bar X^1_{\alpha}-\bar X^2_{\alpha})B^1]\bar x^1_{\alpha}=[-I_n+(I_n-Z)B^1]\alpha z={\bf 0}_n
$, and likewise, $
[-I_n+(I_n-\bar X^1_{\alpha}-\bar X^2_{\alpha})B^2]\bar x^2_{\alpha}=[-I_n+(I_n-Z)B^2](1-\alpha)z={\bf 0}_n
$. Evidently, $(\bar x^1_{\alpha},\bar x^2_{\alpha})=\left(\alpha z,(1-\alpha)z\right)$ satisfies the equilibrium equations for all $\alpha\in[0,1]$. 

We now consider the Jacobian at such an equilibrium. Let $\tilde B_{\alpha}^i$ denote ${\rm{diag}}(B^i\bar x_{\alpha}^i)$.  The Jacobian at any equilibrium on this interval is 
\begin{equation}
J=\begin{bmatrix}
-I_n+(I_n-\bar X^1_{\alpha}-\bar X^2_{\alpha})B^1-\tilde B_{\alpha}^1&-\tilde B_{\alpha}^1\\
-\tilde B_{\alpha}^2&-I_n+(I_n-\bar X^1_{\alpha}-\bar X^2_{\alpha})B^2-\tilde B_{\alpha}^2
\end{bmatrix}
\end{equation}
Now the equilibrium equations give immediately $    {\rm{diag}}(B^1\bar x^1_{\alpha}) =\alpha(I_n-Z)^{-1}Z$ and    ${\rm{diag}}(B^2\bar x^2_{\alpha})=(1-\alpha)(I_n-Z)^{-1}Z$.
% \begin{align*}
%     {\rm{diag}}(B^1\bar x^1_{\alpha})& =\alpha(I_n-Z)^{-1}Z,\\
%     {\rm{diag}}(B^2\bar x^2_{\alpha})& =(1-\alpha)(I_n-Z)^{-1}Z
% \end{align*}
% \[
% {\rm{diag}}(B^1\bar x^1_{\alpha})=\alpha(I-Z)^{-1}Z,\;{\rm{diag}}(B^2\bar x^2_{\alpha})=(1-\alpha)(I-Z)^{-1}Z
% \]
Then, $J$ is seen to be similar to
\small
\begin{align*}
\bar J&=
\begin{bmatrix}
-I_n+(I_n-Z)B^1-\alpha(I_n-Z)^{-1}Z&\alpha(I_n-Z)^{-1}Z\\
(1-\alpha)(I_n-Z)^{-1}Z&-I_n+(I_n-Z)B^2-(1-\alpha)(I_n-Z)^{-1}Z
\end{bmatrix}
\end{align*}
\normalsize
using the transformation matrix $P = {\rm{diag}}(I_n, -I_n)$. The matrix $\bar J$ is an irreducible Metzler matrix. Further, $[z^{\top},\;z^{\top}]^\top \gg {\bf 0}_n$ is a nullvector of $\bar J$ for all $\alpha$. It follows that all other eigenvectors of the matrix have negative real parts. 
It also follows that the bivirus equations defining the line of equilibria define a one-dimensional center manifold along which the Jacobian is singular. By standard center manifold theory, see e.g. \cite[Section 7.6]{sastry1999nonlinearbook}, the eigenvalue properties of $J$ then imply the center manifold is exponentially attractive, i.e. for initial conditions sufficiently close to the manifold, convergence occurs to some point on the manifold exponentially fast. 

\subsection{Proof of \Cref*{thm:boundary_stability}}\label{app:pf_boundary_stability}
We establish the result for the boundary equilibrium $(\bar x^1, {\bf 0}_n)$, with the proof for $({\bf 0}_n, \bar x^2)$ being identical after adjustment of certain indices. First, we give some relevant results concerning the single virus model \eqref{eq:sis}. Let $\mathcal R > 1$, and thus the unique endemic equilibrium $\bar x$ of \eqref{eq:sis} satisfies
\begin{equation}\label{eq:sis_equib}
    (-D + (I_n-\bar X)B)\bar x = {\bf 0}_n.
\end{equation}
The positive vector $\bar x \gg {\bf 0}_n$ can be seen as the nullvector of the matrix $P = -D + (I_n-\bar X)B$, which is an irreducible Metzler matrix. From the discussions in \cref{sec:prelims}, it follows that $s(P) = 0$, which in turn implies that $-P$ is a singular irreducible $M$-matrix. It is known that an irreducible singular $M$-matrix plus a nonnegative diagonal matrix with at least one positive diagonal element yields an irreducible nonsingular $M$-matrix~\cite[Theorem~4.31]{qu2009cooperative_book}. With $\tilde B = \text{diag}(B\bar x)$, it follows that $$-P + \tilde B = D - (I_n-\bar X)B + \tilde B$$ is a nonsingular $M$-matrix, and thus $s(-P+\tilde B) < 0$.

The Jacobian $J(\bar x^1, {\bf 0}_n)$ is upper block triangular, as seen from \eqref{eq:jacobian}. From the arguments immediately above, one can conclude that the upper diagonal block matrix $-I_n + (I_n-\bar X^1)B^1 - \tilde B^1$ is the negative of a nonsingular $M$-matrix, and is therefore Hurwitz. The lower diagonal block matrix $I_n + (I_n-\bar X^1)B^2$ is an irreducible Metzler matrix, and note that $\bar X^1$ is uniquely determined by $B^1$ and is therefore independent of $B^2$. As outlined in \cref{sec:prelims}, one has that $s\big(-I_n + (I_n-Z^1)B^2\big) < 0 \Leftrightarrow \rho\big((I_n-Z^1)B^2\big) < 1$. The condition for instability of $(\bar x^1, {\bf 0}_n)$ can be similarly proved. This completes the proof.

\subsection{Proof of \Cref*{cor:boundary}}\label{app:pf_boundary}

\textit{Item 1:} As detailed below \eqref{eq:sis_equib}, one has that $s(-I_n + (I_n-\bar X^1)B^1) = 0$, which in turn implies that $\rho((I_n-\bar X^1)B^1) = 1$, according to \cref{sec:prelims}. Since $B^2 > B^1$, it follows that there exists at least one entry of $(I_n-\bar X^1)B^2$ is strictly greater than the corresponding entry of $(I_n-\bar X^1)B^1$. Then, \cite[Theorem 2.7]{varga2009matrix_book} establishes that $\rho((I_n-\bar X^1)B^2) > \rho((I_n-\bar X^1)B^1) = 1$,
which in conjunction with Theorem~\ref{thm:boundary_stability} delivers the claim on instability. The argument for local stability of $({\bf 0}_n, \bar x^2)$ is similar, starting with the observation $s(-I_n + (I_n-\bar X^2)B^2) = 0$.

To show there is no equilibrium of the form $(\tilde x^1,\tilde x^2)$ with $\tilde x^1\gg {\bf 0}_n, \tilde x^2\gg{\bf 0}_n$, let us assume to the contrary that such a $(\tilde x^1,\tilde x^2)$ exists. Then, \eqref{eq:equil_1} indicates that $\tilde x^1$ is a positive eigenvector associated with the simple eigenvalue 0 of the irreducible Metzler matrix $-I+(I-\tilde X^1-\tilde X^2)B^1$. Let $y^{\top}$ be the associated positive left eigenvector, normalized to satisfy $y^\top\tilde x^1 = 1$. Set $C=(I-\tilde X^1-\tilde X^2)(B^2-B^1)$. Observe that $C>0_{n\times n}$. Further, \eqref{eq:equil_2} yields $[-I+(I-\tilde X^1-\tilde X^2)B^1]\tilde x^2+C\tilde x^2={\bf 0}_n$.
Premultiplying the left by $y^{\top}$ gives a contradiction, since $C>0_{n\times n},y\gg{\bf{0}}_n, \tilde x^2\gg{\bf{0}}_n$.

\textit{Item 2:} We first show there are no coexistence equilibria, and then show the instability of $(\bar x^1, {\bf 0}_n)$. Suppose, to obtain a contradiction, that there is an interior equilibrium $(\tilde x^1,\tilde x^2)$. Let $J,K$ be indices such that $\tilde x^1_i\leq \tilde x^1_J$ and $\tilde x^2_i\geq \tilde x^2_K$, for all $i$. Observe firstly that $0<1-\tilde x^1_J-\tilde x^2_J\leq 1-\tilde x^1_J-\tilde x^2_K\leq 1-\tilde x^1_K-\tilde x^2_K<1$.
% \begin{equation}\label{eq:useful}
%     0<1-\tilde x^1_J-\tilde x^2_J\leq 1-\tilde x^1_J-\tilde x^2_K\leq 1-\tilde x^1_K-\tilde x^2_K<1.
% \end{equation}
Next, obtain $(1-\tilde x^1_J- \tilde x^2_J)^{-1}\tilde x^1_J=(B^1\tilde x^1)_J\leq (B^1{\bf{1}}_n\tilde x^1_J)_J\leq\bar{b}_1\tilde x^1_J$ from the $J$-th row of \eqref{eq:equil}.
% that
% \begin{align*}
%     (1-\bar x^1_J- \bar x^2_J)^{-1}\bar x^1_J=(B^1\bar x^1)_J\leq (B^1{\bf{1}}_n\bar x^1_J)_J\leq\gamma_1\bar x^1_J
% \end{align*}
It follows that $(1-\tilde x^1_J-\tilde x^2_K)^{-1}\leq (1-\tilde x^1_J-\tilde x^2_J)^{-1}\leq\bar{b}_1$. Similar arguments give $(1-\tilde x^1_J-\tilde x^2_K)^{-1}\geq (1-\tilde x^1_K-\tilde x^2_K)^{-1}\geq \underline{b}_2$ and then there is an immediate contradiction of the hypothesis condition that $\bar{b}_1<\underline{b}_2$.

From $(I-\bar X^1)^{-1}\bar x^1=B^1\bar x^1$ and setting $\bar x^1_J=\max_j\bar x^1_j \in (0,1)$, notice that
$(1-\bar x^1_J)^{-1}\bar x^1_J=(B^1\bar x^1)_J\leq(B^1{\bf{1}}_n \bar x^1_J)_J\leq\bar {b}_1\bar x^1_J<\underline{b}_2\bar x^1_J$
 or equivalently $(1-\bar x^1_J)\underline{b}_2>1$.
%  \begin{equation}\label{eq:yippee}
% (1-\bar x^1_J)\min_i\sum_{j=1}^n\beta^2_{ij}>1
%  \end{equation}
 
For any irreducible nonnegative matrices $A^1,A^2$ the condition $A^1>A^2$ implies $\rho(A_1)>\rho(A_2)$ and $\rho(A^1)\geq \min_i\sum_{j=1}^na_{ij}$. This means that
$\rho[(I-\bar X^1)B^2]\geq\rho[(1-\bar x^1_J)B^2]\geq(1-\bar x^1_J)\underline{b}_2 > 1$
and the instability claim is complete. The proof that $({\bf 0}_n, \bar x^2)$ is locally stable, by establishing $\rho[(I-\bar X^2)B^1]<1$, is the same, \textit{mutatis mutandis}.

\textit{Item 3:} Since $\bar x^2 > \bar x^1$, it immediately follows that $(I_n-\bar X^1)B^2 > (I_n-\bar X^2)B^2$, which according to \cite[Theorem 2.7]{varga2009matrix_book} implies that $\rho((I_n-\bar X^1)B^2) > \rho((I_n-\bar X^2)B^2) = 1$, thus establishing the instability of $(\bar x^1, {\bf 0}_n)$ by the result of Theorem~\ref{thm:boundary_stability}. Similarly, we have that $(I_n-\bar X^1)B^1 > (I_n-\bar X^2)B^1$, and it follows from \cite[Theorem 2.7]{varga2009matrix_book} that $\rho((I_n-\bar X^2)B^1) < \rho((I_n-\bar X^1)B^1) = 1$, and hence $({\bf 0}_n, \bar x^2)$ is locally stable.

\subsection{Proof of \Cref*{thm:trajectoryordering}}\label{app:pf_trajorder}
 
It is evident that with $m=({\bf 0}_n^\top,{\bf 1}_n^\top)$  the initial conditions specified in the theorem statement obey $x_B(0)>_{K_m}x_C(0)>_{K_m}x_A(0)$, and so by \Cref{lem:kamke} on monotone systems, $\phi_t(x_B(0))\gg_{K_m}\phi_t(x_C(0))\gg_{K_m}\phi_t(x_A((0))$ holds for all $t$.
% \begin{equation}
% \phi_t(x_B(0))\gg_{K_m}\phi_t(x_C(0))\gg_{K_m}\phi_t(x_A((0))
% \end{equation}
These inequalities yield the first two lines of \eqref{eq:3trajectory}.  The inequality of the last line was established in \Cref{lem:strictinvariant}. 

{\color{black}\subsection{Proof of \Cref*{lem:global_stab}}\label{app:pf_global_stab}
%  { \color{red} If the cyan inserts are accepted, this proof comes out. }
 Without loss of generality, suppose that $\mathcal{S}$ has a globally stable equilibrium at $\bar x$. To obtain a contradiction, suppose there is an initial state $x(0)\in\Delta$ for $\hat{\mathcal{S}}$ for which the trajectory $x(t)$ does not converge to $\bar x$. Let $\mathcal N\subset \Delta$ be a suitably small neighborhood of $x(0)$. Then there exist $x_A(0)\in\mathcal N$ and $x_B(0)\in\mathcal N$ with $x_B(0)\gg x(0)\gg x_A(0)$ such that, by \Cref{thm:convergence}, the trajectories $x_A(t)$ and $x_B(t)$ starting from $x_A(0)$ and $x_B(0)$ converge to a stable equilibrium of $\hat{\mathcal{S}}$.  But by \Cref{prop:Didentity}, the system $\hat{\mathcal{S}}$ cannot have a stable equilibrium other than $\bar x$, since the system $\mathcal{S}$ has no stable equilibrium point other than $\bar x$. In other words, convergence of trajectories in the system $\hat{\mathcal{S}}$ must occur to $\bar x$. Now use the key result of \Cref{thm:trajectoryordering} on $\hat{\mathcal{S}}$. The inequality for the initial conditions implies $x_B(t)\gg x(t)\gg x_A(t)\;\forall t$ and taking the limit as $t\to\infty$ yields $\lim_{t\to\infty}x(t)=\bar x$. Note we can always select $x_A(0)$ and $x_B(0)$ to satisfy \Cref{thm:trajectoryordering}; due to \Cref{lem:strictinvariant}, any $x(0) \in \Delta$ will satisfy ${\bf 0}_n \ll x^i(t) \ll {\bf 1}_n$, $i = 1,2$, after a finite time $t$.}

\subsection{Proof of \Cref*{cor:equilibriumordering}}\label{app:pf_equiborder}

 The underlying system is an irreducible monotone system and by assumption, almost all trajectories converge to a limit point  (which may be locally stable, or a saddle point); the remaining trajectories are nonattractive limit cycles. The inequalities established in \Cref{thm:trajectoryordering} imply that the limiting trajectory resulting from $x_C(0)$ necessarily lies within the hyperrectangle $\mathcal W$, and this can only be a nonattractive limit cycle if $x_C(0)$ itself lies in $\mathcal W$ and $\mathcal W$ is neither a single point nor a one-dimensional interval. These observations imply one of the three conclusions in the corollary statement must hold. 
 
\subsection{Proof of \cref*{cor:no_interior_equib}}\label{app:pf_no_interior_equib}
The argument relies on the fact that all trajectories of the bivirus system approach a locally stable equilibrium, except possibly from a set of initial conditions of measure zero, as per \cref{thm:convergence} and \cref{lem:monotone_convergence}. If there were two unstable boundary equilibria, a contradiction is immediate because almost all trajectories approach a stable equilibrium; boundary equilibria are excluded through their instability and a coexistence equilibrium is not present by hypothesis. If there were two stable boundary equilibria, \cref{cor:equilibriumordering} yields existence of an interior unstable equilibrium, another contradiction. Hence there is precisely one stable boundary equilibrium. Suppose that there exists an interior point $x(0) = (x^1(0),x^2(0))$ for which the associated trajectory does not converge to the stable boundary equilibrium (perhaps on a nonattractive limit cycle). Let $\mathcal N$ be a suitably small neighbourhood of $x(0)$. Then, there exist $x_A(0)\in\mathcal N$ and $x_B(0)\in\mathcal N$ with $x_B(0)\gg x(0)\gg x_A(0)$ such that the trajectories starting from $x_A(0)$ and $x_B(0)$ converge to the stable boundary equilibrium, and \cref{cor:equilibriumordering} implies that the trajectory beginning at $x(0)$ also converges to the common limit, delivering the contradiction. 
 
\section{Applying algebraic geometry to the proof of \Cref*{thm:convergence}}
The main focus is to prove that for generic $D^i, B^i$, the bivirus system in \eqref{eq:twovirus} has a finite number of equilibria. To do this, we shall first argue that there can only be an infinite number of equilibria for values of the free parameters in $D^i,B^i$ lying in an algebraic set, i.e. a set defined by setting a multivariate polynomial in the parameters to zero. The existence of this set will be demonstrated by algebraic geometry. If the multivariate polynomial is the zero polynomial, the algebraic set comprises the whole space. If it is not zero, then the algebraic set has measure zero and for almost all parameter values, there will be a finite number of equilibria.  Then we prove that for any given $n$, there exists a specific choice of $D^i$ and $B^i$ for which the bivirus system has a finite number of equilibria, implying that the algebraic set cannot be the whole set and so necessarily has measure zero. 
 
\subsection{Background on algebraic geometry}
We first provide relevant background on algebraic geometry, the main reference being \cite[Chapter 3]{cox2006using}. 
% The key algebraic geometry result is a generalization of some basic ideas of linear algebra, which we first review. Recall that the linear equation $Ax=b$ with $A\in\mathbb R^{n\times n}, b\in\mathbb R^n$ and $x\in\mathbb R^n$ the unknown to be found has a unique solution provided that a certain expression polynomial{\footnote{ The expression is in fact multilinear, but this is not of concern to us.}} in the entries of $A$ (the determinant in fact) is nonzero, and otherwise has no solutions or an infinity of solutions. Sometimes, $A$ has some structure in the sense that certain entries assume fixed values, and other entries contain free or unknown parameters, a companion matrix being an example. (Free parameters may not always be completely free, for example, they might be constrained to lie in an interval.)  In this case, it is often convenient to regard the determinant as a multivariate polynomial in the free parameters. If $\beta$ is the vector of free parameters, the condition for unique solvability is of the form $\Delta(\beta)\neq 0$ where $\Delta$ is the determinant of $A$. 

%  Algebraic geometry generalizes this result. 
With $x$ denoting an $n$-vector of unknowns, a \textit{monomial} is a term of the form $x_1^{\alpha_1}x_2^{\alpha_2}\dots x_n^{\alpha_n}$, where the $\alpha_i$ are nonnegative integers. A single multivariate polynomial equation is obtained by taking a linear combination of such monomials and setting it to zero. Algebraic geometry enables examination of the solvability question for $n$ such equations in $n$  unknowns. They can be written in the form  $P(x,\beta)=0$ where $\beta$ is the set of free coefficients of the different monomials. (In any one equation, some monomials may be absent, and others again may have a fixed coefficient, such as 1; the coefficients of the remaining monomials make up the entries of $\beta$). We focus on the case of real polynomial equations, so that $\beta$ can be regarded as a real vector, the first block of entries corresponding to monomials in the first equation, the second block to monomials in the second equation, and so on. 
 
 In order to deal with a phenomenon (so-called `solutions at infinity') not encountered in working with a purely linear equation set, a technical device is introduced. 
 
A single multivariate polyonomial equation is termed homogeneous when the sum of the powers in each monomial in the equation is the same for all terms. If one has a nonhomogeneous equation, it can be made homogeneous by introducing a further variable, $x_0$, with a monomial $x_1^{\alpha_1}x_2^{\alpha_2}\dots x_n^{\alpha_n}$ being replaced by $x_0^{\alpha_0}x_1^{\alpha_1}x_2^{\alpha_2}\dots x_n^{\alpha_n}$ with $\alpha_0$ chosen to ensure all monomials in the single polynominal equation have the same total degree $\sum_{i=0}^n\alpha_i$. 
 
 Consider in fact $n$ such homogeneous polynomial equations in the $n+1$ unknowns $x_0,x_1,\dots x_n$. Note that while separately homogeneous, the individual equation degrees are not necessarily the same. The equations may well have been obtained by making homogeneous the inhomogeneous set $P(x,\beta)=0$. We shall write the homogeneous set as $\bar P(x_0,x,\beta)=0$. This is related to $P$ by $\bar P(1,x,\beta)=P(x,\beta)$.  The key result from algebraic geometry is the following, \cite[pg. 86, Theorem 2.3]{cox2006using}. 
 
 \begin{theorem}\label{thm:B1}
 Consider a set of $n$ homogeneous polynomial equations, denoted $\bar P(x_0,x,\beta)=0$ with free parameters in the vector $\beta$ and with unknowns $x_0,x_1,\dots,x_n$. Then there exists an expression $R(\beta)$ (termed the resultant) that is polynomial in the entries of $\beta$ such that if $R(\beta) \neq 0$ for particular values of the $\beta$ entries, there are a finite number of nonzero solutions (possibly complex) to $\bar P(x_0,x,\beta)=0$, disregarding scaling. If $R(\beta) = 0$, there are either no solutions or an infinity of solutions (disregarding scaling).
 \end{theorem}
 
 Evidently, if $(\bar x_0,\bar x_1,\dots,\bar x_n)$ is a solution to a homogeneous set of equations $\bar P(x_0,x,\beta)=0$, then $(\lambda\bar x_0,\lambda\bar x_1,\dots\lambda\bar x_n)$ for any nonzero real or complex $\lambda$ is also a solution. If $\bar x_0\neq 0$, the choice $\lambda=(\bar x_0)^{-1}$ recovers a solution to the inhomogeneous set $P(x,\beta)=0$. However, if $\bar x_0=0$, a solution to $P(x,\beta)=0$ cannot be recovered in this way; such solutions are however sometimes termed `solutions at infinity' of the set $P(x,\beta)=0$, \cite[pg. 115]{cox2006using}.
 
\begin{corollary}
Assume the same hypotheses as Theorem \ref{thm:B1}.
If for a particular choice of values for the entries of $\beta$, call it $\hat\beta$, the resultant takes a  nonzero value, then for almost all choices of $\beta$, the resultant will be nonzero.
\end{corollary}

The proof is immediate: it is obvious that if a polynomial in a single variable is nonzero for some value of that variable, it is nonzero for almost all values, i.e.
nonzero everywhere except on a set of measure zero (which in the case of a polynomial in a single variable is a finite set). The first property is clearly also true for any (nonzero) multivariate polynomial such as the resultant.

Thus if the resultant takes a nonzero value for a particular choice of values $\hat \beta$ for the free parameters of the equation set $\bar P(x_0,x,\beta)$, then for almost all values assumed by the entries of $\beta$, the resultant polynomial $R(\beta)$  will evaluate as a nonzero number. Further, if it takes a nonzero value for a particular choice of values for the free parameters in  the equation set $\bar P(x_0,x,\beta)=0$ and there is an associated solution with $x_0\neq 0$, i.e. there is a solution of $P(x,\beta)=0$, the same will hold true for almost all values of the free parameters. 

%  Given a set of multivariate polynomials with associated free parameter set defined by the entries of a vector $\beta$, a property is said to \textit{hold generically} \cite[see page 115]{cox2006using} if there is a nonzero polynomial in the entries of $\beta$ such that the property holds for all particular values of $\beta$ for which the polynomial is nonzero. Evidently then, the equation set $P(x,\beta) =0$, obtained by setting the polynomials to zero, generically has a finite number of solutions if the resultant $R(\beta)$ is nonzero. 
 
 \subsection{Proof of \Cref*{thm:convergence}}\label{app:pf_convergence}
As noted just prior to \cref{def:free_param}, evidently the key part of the theorem we have to prove is that the number of equilibria is finite.

The equilibrium equations for the bivirus system, given in \eqref{eq:equil}, are a set of polynomial equations in the $\bar x_i^j$. Suppose that we take $D^1=D^2=I_n$, and $B^1,B^2$ as positive diagonal matrices with no two corresponding entries equal. The equation set  then becomes decoupled, and the $i$-th entry of each of the vector equations in \eqref{eq:equil} is
\begin{align*}
    [-1+(1-\bar x^1_i-\bar x^2_i)\beta^1_{ii}]\bar x^1_i&=0\\
    [-1+(1-\bar x^1_i-\bar x^2_i)\beta^2_{ii}]\bar x^2_i&=0
\end{align*}
The solvability of these equations yielding a finite number of solution pairs $(\bar x^1_i,\bar x^2_i)$ is easily checked. The solutions in fact are $( 0,0),(0,1-(1/\beta^2_{ii})),(1-(1/\beta^1_{ii}),0)$. Of course, the other entries of the $\bar x^i$ can be treated in the same way due to the decoupling. Since there is a particular choice for the entries of $B^1,B^2$ and $D^1, D^2$ (the latter particular choice being the identity matrix) yielding a finite number of solutions, the algebraic geometry arguments presented above show that for almost all choices of the free parameters in $B^i$, and the diagonal entries of $D^i$, a finite number of solutions to the equilibrium equations exist.

{\color{black}Due to \Cref{lem:strictinvariant}, it follows that for any $x(0) \in \Delta$, we have $x(\tau) \in \tilde \Delta$ for some finite time $\tau \geq 0$, with $\tilde \Delta$ an open strict subset of $\Delta$ that is positively invariant. With finiteness of the equilibria assured, we apply \Cref{lem:monotone_convergence} to the bivirus system \eqref{eq:twovirus}, relating $\mathcal{M}$ to $\tilde \Delta$. Note, we can relate $\overline{\mathcal M}$ to $\overline{\Delta} \triangleq \{x^1, x^2 \vert {\bf 0}_n \leq x^i \leq {\bf 1}_n  \text{ for } i = 1,2, \text{and } x^1+x^2 \leq {\bf 1}_n\}$. }
% \Cref{lem:strictinvariant} establishes trajectories starting on the boundary of $\overline{\Delta}\setminus $ that $\overline{\Delta} \setminus (\Delta\cup {\bf 0}_n) $ is \textit{not} a positive 
% Since there is a particular choice for the entries of $B^1,B^2$, even with the specialized choice of $D^1,D^2$ as identity matrices, yielding a finite number of solutions,  the algebraic geometry arguments presented above show that for almost all choices of the free parameters in $B^i$, and, if desired, the diagonal entries of $D^i$, a finite number of solutions to the equilibrium equations exists. 

\section*{Acknowledgments}
{\color{black}We would like to acknowledge the anonymous referees and their comments and suggestions for improving the manuscript.}

\bibliographystyle{siamplain}
\bibliography{MYE_ANU,ji}
\end{document}